\documentclass[12pt]{article}
\usepackage{mathrsfs}
\usepackage{amsfonts}
\usepackage{psfrag,epsfig}
\usepackage{latexsym,amsmath,amssymb}
\usepackage{exscale,relsize}
\usepackage{siunitx}
\usepackage{graphicx}
\usepackage{epstopdf}
\usepackage{indentfirst,cite}
\numberwithin{equation}{section}

\begin{document}
\title{{Existence of positive solutions for Kirchhoff type problems with critical exponent in exterior domains}
\thanks{Research supported by National Natural Science Foundation of China (No. 12001403, 11571187, 11771182). The corresponding author is Shiwang Ma.}
\author {{Liqian Jia$^a$, \quad Xinfu Li$^b$,\quad Shiwang Ma$^c$}\\
\small \it $^a$School of Mathematics and Statistics, Weifang University, Weifang 261061, P. R. China\\
\small \it $^b$School of Science, Tianjin University of Commerce,
Tianjin 300134, P. R. China\\
\small \it $^c$School of Mathematical Science and LPMC, Nankai University, Tianjin 300071, P. R. China\\
\small E-mails: jialiqian1992@163.com, lxylxf@tjcu.edu.cn,
shiwangm@nankai.edu.cn.\ \ }}

\date{}
\maketitle \baselineskip=12pt \noindent {\small {\bf Abstract}: 
In this paper, by using variational methods  we study  the existence of positive solutions for the following Kirchhoff type problem:
$$
\left\{
\begin{array}{ll}
-\left(a+b\mathlarger{\int}_{\Omega}|\nabla u|^{2}dx\right)\Delta u+V(x)u=u^{5}, \ & x\in\Omega,\\
\\
u=0,\ & x\in\partial \Omega,
\end{array}\right.
$$
where $a>0$, $b\geq0$, $\Omega\subset\mathbb R^3$  is an unbounded exterior domain, $\partial\Omega\neq\emptyset$, $\mathbb{R}^{3}\backslash\Omega$ is bounded,  $u\in D_{0}^{1,2}(\Omega)$, and $V\in L^{\frac{3}{2}}(\Omega)$ is a  non-negative continuous function.  It turns out  that the above Kirchhoff equation has no ground state solution. Nonetheless,  by establishing some global compact lemma and constructing a suitable minimax value $c$ at a higher energy level where so called Palais-Smale condition holds, we succeed to obtain a positive solution for such a problem whenever  $V$ and the hole $\mathbb{R}^{3}\setminus\Omega$ are suitable small in some senses.
To the best of our knowledge, there are few similar results published in the literature concerning the existence of positive solutions for  Kirchhoff equation in exterior domains. Our result also holds true in the case $\Omega=\mathbb R^3$, particularly, if $a=1$ and $b=0$,  we  improve  some existing results (such as Benci, Cerami, Existence of positive solutions of the equation $-\Delta u+a(x)u=u^{(N+2)/(N-2)}$ in $\emph{R}^{N}$, J. Funct. Anal., 88 (1990), 90--117) for the corresponding Schr\"odinger equation in the whole space.

\vskip0.5cm \noindent {\bf AMS (2020)  Subject Classification}:
35J60; 35J20; 35A01; 35B09; 35B33.

\baselineskip=20pt

\section{Introduction and main results}
\setcounter{equation}{0}

In this paper, we consider the existence of  positive solutions to the following
Kirchhoff type problem:
\begin{equation}\label{eq:Omega}
  \left\{
\begin{array}{ll}
-\left(a+b\mathlarger{\int}_{\Omega}|\nabla u|^{2}dx\right)\Delta u+V(x)u=f(u), &\quad x\in\Omega,\\
\\
u=0, &\quad x\in\partial \Omega,\\
\end{array}\right.
\end{equation}
where $a>0$, $b\geq0$, $\Omega\subset\mathbb R^3$  is
an unbounded exterior domain, $\partial\Omega\neq\emptyset$,
$\mathbb R^3\backslash\Omega$ is bounded, $u\in D_{0}^{1,2}(\Omega)$, $f(u)=u^5$, 
and $V\in L^{\frac{3}{2}}(\Omega)$ is a given non-negative continuous function.

The Kirchhoff type problems arose in the following time dependent wave
equation proposed by Kirchhoff (see \cite{Kirchhoff})
$$
\rho\frac{\partial^{2}u}{\partial t^{2}}-\Big(\frac{P_{0}}{h}+\frac{E}{2L}\int_{0}^{L}\Big|\frac{\partial u}{\partial x}\Big|^{2}\Big)\frac{\partial^{2}u}{\partial x^{2}}=0,
$$
where $L$ is the length of the string, $h$ is the area of the cross
section, $E$ is the Young modulus of the material, $\rho>0$ is the
mass density and $P_{0}$ is the initial tension. It is an extension
of the classical d'Alembert wave equation for free vibration of
elastic strings. In a classical paper \cite{Lions}, J. L. Lions introduced an abstract functional
framework to the equation
$$
\left\{
\begin{array}{ll}
u_{tt}-(a+b\int_{\Theta}|\nabla u|^{2}dx)\Delta u=f(x,u), \quad x\in \Theta,\\
\\
u=0,\ \ \quad\qquad\qquad\qquad\qquad\qquad\qquad\quad x\in\partial \Theta,\\
\end{array}\right.
$$
where $\Theta$ is a bounded domain of $\mathbb{R}^{3}$. After that,
people pay much attention to the Kirchhoff problems,
see \cite{Arosio,Cavalcanti,D'Ancona} and the references therein.

Another interest in studying the Kirchhoff type problem comes from the non-locality which means that the
problem is no longer a pointwise identity. This fact make  the study of the Kirchhoff problem is quite different from
the local elliptic problems.  
In recent years, the Kirchhoff type problems defined in various bounded domains
have been studied
extensively,  and we refer
the readers to
\cite{Alves-1,Bensedki,Chen,He-1,Mao,Yang,Zhang}  and references therein.
The Kirchhoff problem defined on the whole space
 \begin{equation}\label{eq:f(x,u)-RN}
-\left(a+b\mathlarger{\int}_{\mathbb{R}^{N}}|\nabla u|^{2}dx\right)\Delta u+V(x)u=f(u), \quad \text{in} \ \ \mathbb{R}^{N}
\end{equation} 
has been studied by many authors when $V(x)$ and $f(u)$ satisfy various suitable conditions. 
For a subcritical nonlinearity $f(u)$, 
we refer the readers to \cite{Guo,He-3,Li-0,Li-1,Liu-0,Liu-3,Lu,Sun-1,Sun-2} for existence results of solutions to \eqref{eq:f(x,u)-RN} and its variants. 
Many existence results  for \eqref{eq:f(x,u)-RN}  with a critical nonlinearity $f(u)$ have also been obtained, see
\cite{He-4,He-5,He-6,Li-Xiang,Xie-1,Liu-1,Liu-2,Liu-Luo,Wang} and the references therein.  Particularly, if $N=3$, $f(u)=u^5$ and $V (x) \ge  0, \forall x\in \mathbb R^3$,  Xie et al. \cite{Xie-1} obtained a bound state solution to \eqref{eq:f(x,u)-RN} under  the condition 
\begin{equation}\label{e1.1}
0<|V|_{\frac{3}{2}}<V_{a,b,S},\ \ V_{a,b,S}:=\left(2aK^{4}+\frac{b}{2}K^{5}S^{\frac{3}{2}}\right)^{\frac{1}{3}}S-K^{2}S,
\end{equation}
where and in what follows, $K:=\frac{1}{2}\left(bS^{\frac{3}{2}}+\sqrt{b^{2}S^{3}+4a}\right)>0$, $S$ is the best Sobolev constant for the embedding $D^{1,2}(\emph{R}^{3})\hookrightarrow L^{6}(\emph{R}^{3})$, that is,
\begin{equation}\label{e1.11}
S:=\inf_{u\in D^{1,2}(\mathbb{R}^3)\backslash \{0\}}\frac{\int_{\mathbb{R}^3}|\nabla u|^2dx}{\left(\int_{\mathbb{R}^3}u^6dx\right)^{1/3}}.
\end{equation}
 
 However, when the domain $\Omega$ is an unbounded exterior domain,  the existence of solutions to elliptic  problems become  more delicate and difficult. 
 There are  a few papers in the  literature dealing with  local elliptic problems in exterior domain.  The authors in \cite{Alves2017,Li2004,Montenegro,Hirano} studied the following Schr\"odinger equation
\begin{equation}\label{eq:elliptic}
-\Delta u+V(x)u=f(u), \quad x\in \Omega\subset\mathbb R^N,\quad N\geq3
\end{equation}
under various conditions. Specifically,  Alves and de Freitas \cite{Alves2017} considered \eqref{eq:elliptic} in the case $N=3$, $V(x)=1$, $f(u)=g(u)+\varepsilon u^{5}$ with $g(u)$ being subcritical and $\varepsilon>0$ being small.
Li et al. \cite{Li2004} and Hirano \cite{Hirano} studied \eqref{eq:elliptic} with $N\geq 3$, $V(x)=0$, $f(u)=u^{2^{*}-1}$, where $2^{*}=2N/(N-2)$ is the Sobolev critical exponent. When $N\geq3$, $V(x)=1$, $f(u)=g(u)+u^{2^{*}-1}$,  Montenegro and  Abreu \cite{Montenegro} studied \eqref{eq:elliptic} with Neumann boundary condition in a radially symmetric exterior domain.
 
 To the best of our knowledge, a few authors consider the  Kirchhoff type equations in exterior domains and little is achieved until now \cite{Chen-Liu,Jia,Wang-1}.
 When $\Omega=\mathbb R^3$ and $V(x)\equiv \lambda$,   the scaling
$$
v(x)=u(\sqrt{\varpi_u}x), \quad \varpi_u=a+b\int_{\mathbb R^3}|\nabla u|^2dx,
$$
reduces the nonlocal problem \eqref{eq:Omega} to the local problem
\begin{equation}\label{eq:local}
-\Delta v+\lambda v=f(v), \quad x\in \mathbb R^3.
\end{equation}
However, under our assumption $\Omega\not=\mathbb R^3$, $V(x)\not\equiv {\rm const}$, the situation is  quite different and above reduction does not work any more. 
 Recently, authors in  \cite{Chen-Liu,Jia} considered the Kirchhoff problem \eqref{eq:Omega} in exterior domains
under the assumptions  that $V(x)\equiv\lambda$  is a positive constant and $f$ satisfies $\lim_{|s|\rightarrow+\infty}f(s)/s^{5}=0$. 
Wang et al. \cite{Wang-1} extended the result in \cite{Alves2017} to Kirchhoff problem  and obtained a positive solution to \eqref{eq:Omega}  when $V(x)\equiv 1$ and $f(u)=u^{q-1}+\varepsilon u^5$, with $q\in (4,6)$ and $\varepsilon>0$ being sufficiently small.

In the present paper, we assume that the non-negetive function $V(x)\not\equiv {\rm const}$, $\Omega\not=\mathbb R^3$ and $f(u)=u^5$, that is, $f(u)$ is a critical nonlinearity.  Both the unboundedness of exterior domain and the critical nonlinearity simultaneously make the problem lack of compactness, so the standard variational methods and techniques do not apply directly.  Inspired by the results in \cite{Benci-1, Xie-1}, by establishing some global compact lemma and constructing a suitable minimax value $c$ at a higher energy level where so called Palais-Smale condition holds, we succeed to obtain a positive solution for such a problem whenever  $V$ and the hole $\mathbb{R}^{3}\setminus\Omega$ are suitable small in some senses. To state our main result, we make the following assumption. 

\smallskip
\noindent$\bf{(V_{1})}\quad$ $V\in L^{\frac{3}{2}}(\Omega)$ is a non-negative continuous function and satisfies \begin{equation}\label{con:V}
0<|V|_{\frac{3}{2}}<C_{a,b,S},
\end{equation}
where 
$$
C_{a,b,S}:=\left(2aKS^{\frac{3}{2}}+\frac{b}{2}K^{2}S^{3}\right)^{\frac{2}{3}}-bS^{2}\left(2aKS^{\frac{3}{2}}+\frac{b}{2}K^{2}S^{3}\right)^{\frac{1}{3}}-aS.
$$

\smallskip 

 Our main result is as follows:\\
\textbf{Theorem 1.1.} {\it Let $a>0$, $b\geq0$, $f(u)=u^5$ and $(V_{1})$ hold.
Then there exists a small $\hat{\rho}>0$ such that
if
$$
{\rm diam}(\mathbb{R}^{3}\setminus\Omega):=\sup\left\{|x-y|: \ x,y\in \mathbb{R}^{3}\setminus\Omega\right\}<\hat{\rho},
$$
the problem \eqref{eq:Omega} admits at least one positive solution.}

\smallskip

\noindent{\bf Remark 1.1.} In this paper, we adapt nontrivially the arguments used in \cite{Benci-1, Xie-1} because of the presence of the nonlocal term and the nontrivial domain $\Omega$, which make the problem  more complicated and involved.  
It is easy to see that Theorem 1.1 is still valid for $\Omega=\mathbb R^3$ and may be viewed as generalization of the main results in \cite{Benci-1, Xie-1}.  Particularly, we greatly  improve  the
main result in \cite{Benci-1} in the case $N = 3$. 
In fact, if $a = 1$, $b = 0$ and $\Omega=\mathbb R^3$, then \eqref{eq:Omega} reduces to
\begin{equation}\label{eq:local2}
-\Delta u + V (x)u = u^5, \  \  in\  \ \mathbb R^3, 
\end{equation}
and \eqref{con:V} reduces to 
\begin{equation}\label{eq:local3}
0 < |V |_{\frac{3}{2}}< S(2^{2/3} -1).
\end{equation}
 In the excellent paper \cite{Benci-1}, Benci and Cerami  got a positive solution of \eqref{eq:local2}
in $D^{1,2}(\mathbb R^3)$ under  \eqref{eq:local3} and the following additional assumptions: $V (x) \ge  0$ for all $x\in \mathbb R^3$, there exist two positive
numbers $p_1 < 3/2 < p_2 < 3$ such that 
\begin{equation}\label{eq:local4}
V\in L^p(\mathbb R^3) \quad      for \   all \ \   p\in [p_1,p_2].
\end{equation}
In this paper, without  the condition \eqref{eq:local4},  we  obtain a positive solution of \eqref{eq:local2} and recover the main result in \cite{Benci-1}. 

\noindent{\bf Remark 1.2.} We also remark that our main result improves the result in \cite{Xie-1} for small $b>0$. In fact, our assumption \eqref{con:V} is weaker than the condition \eqref{e1.1} for small $b>0$.  To see this, we observe that  if  $b>0$ small enough, then for any fixed $a>0$, $K=a^{\frac{1}{2}}+o_b(1)$ and hence 
$$
C_{a,b,S}=(2^{\frac{2}{3}}-1)aS+o_b(1)>0, \quad V_{a,b,S}=(2^{\frac{1}{3}}-1)aS+o_b(1)>0,
$$
where $\lim_{b\to 0}o_b(1)=0$, which yields that $C_{a,b,S}>V_{a,b,S}$ for small $b>0$. 

This paper is organized as follows. In Section 2, we establish the variational framework of \eqref{eq:Omega} and give some preliminary lemmas. In Section 3,
 we give the compactness results which are needed  in the proof of  our  main result. In Section 4, we give the detailed proof of
  Theorem 1.1.

\vskip 5mm

\section{Preliminary lemmas}

In this section, we give the variational setting and some notations for the problem \eqref{eq:Omega}. Without any loss of generality, we may assume that $0\in\emph{R}^{3}\setminus\Omega$.
As usual, $H^{1}(\mathbb{R}^{3})$ is the Sobolev space endowed with the inner product and norm, respectively,
$$
(u,v)_{H^{1}}=\int_{\mathbb{R}^{3}}(\nabla u\nabla v+uv)dx,\quad \|u\|_{H^{1}}^{2}=(u,u)_{H^{1}}.
$$
Let $1\leq q<+\infty$, $L^{q}(\mathbb{R}^{3})$ and $L^{q}(\Omega)$ are the class
of measurable functions $u$ which satisfies $\int_{\mathbb{R}^{3}}|u|^{q}dx<+\infty$ and $\int_{\Omega}|u|^{q}dx<+\infty$
respectively,  and  the norms in $L^{q}(\mathbb{R}^{3})$ and $L^{q}(\Omega)$ are given  by
$$
\big|u\big|_{L^{q}(\mathbb{R}^{3})}=\left(\int_{\mathbb{R}^{3}}|u|^{q}dx\right)^{\frac{1}{q}},\quad |u|_{q}=\left(\int_{\Omega}|u|^{q}dx\right)^{\frac{1}{q}},\quad
$$
respectively.

We consider the equation
\eqref{eq:Omega} in the space $D^{1,2}_{0}(\Omega)$. We recall that for any smooth open domain $D\subset\mathbb{R}^{3}$,
$D^{1,2}_{0}(D)$ denotes the Hilbert space defined by
$$
D^{1,2}_{0}(D):=\left\{u\in L^{6}(D):|\nabla u|\in L^{2}(D)\right\}
$$
with the inner product and norm, respectively,
$$
(u,v):=\int_{D}\nabla u\nabla vdx,\quad \|u\|^{2}=(u,u).
$$
And we denote respectively the inner product and norm of the Hilbert space
$D^{1,2}(\mathbb{R}^{3}):=\left\{u\in L^{6}(\mathbb{R}^{3}):|\nabla u|\in L^{2}(\mathbb{R}^{3})\right\}$ by
$$
(u,v)_{\mathbb{R}^{3}}:=\int_{\mathbb{R}^{3}}\nabla u\nabla vdx,\quad \|u\|_{\mathbb{R}^{3}}^{2}=(u,u)_{\mathbb{R}^{3}}.
$$
It is well known that the solutions of \eqref{eq:Omega} are the critical points of the energy functional $I: D^{1,2}_{0}(\Omega)\rightarrow\mathbb{R}$
defined by
\begin{equation}\label{eq:I(u)}
I(u):=\frac{a}{2}\int_{\Omega}|\nabla u|^{2}dx+\frac{b}{4}\left(\int_{\Omega}|\nabla
u|^{2}dx\right)^{2}+\frac{1}{2}\int_{\Omega}V(x)u^{2}dx-\frac{1}{6}\int_{\Omega}u^{6}dx.
\end{equation}
We define the Nehari manifold $\mathcal{N}$ for the
functional $I$ by
$$
\mathcal{N}:=\left\{u\in D_{0}^{1,2}(\Omega)\backslash\{0\}:G(u)=0\right\}, \quad
G(u):=\langle I'(u),u\rangle,
$$
and we give the following lemma to state some properties of $\mathcal{N}$.
\\
\textbf{Lemma 2.1.} {\it Let $a> 0$, $b\geq0$, $V\in L^{\frac{3}{2}}(\Omega)$ be a non-negative function, then the following statements hold:\\
$(1)$ $\mathcal{N}$ is a $C^{1}$ regular manifold diffeomorphic to the unit sphere of $D^{1,2}_{0}(\Omega)$;\\
$(2)$ $I$ is bounded from below on $\mathcal{N}$ by a positive constant;\\
$(3)$ $u$ is a critical point of $I$ if and only if $u$ is a critical point of $I$ constrained on $\mathcal{N}$.}\\
\textbf{Proof.} $(1)$ Let $u\in D^{1,2}_{0}(\Omega)$ be such that
$\|u\|=1$. We claim that there exists a unique $t\in (0,+\infty)$
such that $tu\in\mathcal{N}$. In fact, let $\alpha_{u}(t):=I(tu)$.  For any $t>0$, we observe that
$$
\alpha_{u}(t)=\frac{1}{2}t^{2}\int_{\Omega}\left(a|\nabla u|^{2}+V(x)u^{2}\right)dx+\frac{1}{4}bt^{4}\left(\int_{\Omega}|\nabla u|^{2}dx\right)^{2}-\frac{1}{6}t^{6}\int_{\Omega}u^{6}dx.
$$
Clearly, $\alpha_{u}(0)=0$ and
$\alpha_{u}'(0)=0$. By $a>0$ and $V(x)\geq 0$, we have $\alpha_{u}(t)>0$ for $t>0$ small enough and $\alpha_{u}(t)<0$
for $t>0$ large enough. Hence, $\max_{t\geq0}\alpha_{u}(t)$ is
achieved at a $t_{u}>0$ such that $\alpha_{u}'(t_{u})=0$, and the
corresponding point $t_{u}u\in\mathcal{N}$ which is called the
projection of $u$ on $\mathcal{N}$. Suppose by contradiction  that
there exist $t_{1}>t_{2}>0$ such that $t_{1}u,
t_{2}u\in\mathcal{N}$, then we have
$$
\frac{1}{t_{1}^{2}}\int_{\Omega}\left(a|\nabla u|^{2}+V(x) u^{2}\right)dx+b\left(\int_{\Omega}|\nabla u|^{2}dx\right)^{2}=t_{1}^{2}\int_{\Omega}u^{6}dx,
$$
$$
\frac{1}{t_{2}^{2}}\int_{\Omega}\left(a|\nabla u|^{2}+V(x) u^{2}\right)dx+b\left(\int_{\Omega}|\nabla u|^{2}dx\right)^{2}=t_{2}^{2}\int_{\Omega}u^{6}dx,
$$
which gives that
$$
0>\left(\frac{1}{t_{1}^{2}}-\frac{1}{t_{2}^{2}}\right)\int_{\Omega}\left(a|\nabla u|^{2}+V(x) u^{2}\right)dx=(t_{1}^{2}-t_{2}^{2})\int_{\Omega}u^{6}dx>0.
$$
That is a contradiction. Thus, $t_{u}$ is unique such that
$t_{u}u\in\mathcal{N}$. Moreover, $\alpha'_{u}(t)>0$ for $0<t<t_{u}$,
$\alpha'_{u}(t)<0$ for $t>t_{u}$.

Moreover, it is well known that $I\in
C^{2}\left(D^{1,2}_{0}(\Omega),\mathbb{R}\right)$, thus $G$ is a $C^{1}$
functional. For any $u\in\mathcal{N}$, we have
\begin{equation}\label{eq:<G'(u),u>}
\begin{array}{rcl}
\langle G'(u),u\rangle&=&\mathlarger{2a\|u\|^{2}+4b\|u\|^{4}+2\int_{\Omega}V(x) u^{2}dx-6\int_{\Omega}u^{6}dx}\\
\\
&=&\mathlarger{2\left(-a\|u\|^{2}-\int_{\Omega}V(x) u^{2}dx-\int_{\Omega}u^{6}dx\right)}<0.
\end{array}
\end{equation}
So (1) is proved.

$(2)$ Let $u\in\mathcal{N}$. By $V(x)\geq 0$, $a>0$, $b\geq 0$ and the Sobolev inequality, there exists  $k_1>0$  such that
$$
0=a\|u\|^{2}+\int_{\Omega}V(x) u^{2}dx+b\|u\|^{4}-\int_{\Omega}u^{6}dx\geq a\|u\|^{2}-k_{1}\|u\|^{6}.
$$
Hence, there exists $k>0$ such that
\begin{equation}\label{eq:u lower bound}
\|u\|\geq k>0.
\end{equation}
This implies that
$$
I|_{\mathcal{N}}(u)=\mathlarger{\frac{a}{3}\|u\|^{2}+\frac{1}{3}\int_{\Omega}V(x) u^{2}dx+\frac{b}{12}\|u\|^{4}}
\geq\mathlarger{\frac{a}{3}k^{2}+\frac{b}{12}k^{4}}>0.
$$
Therefore, (2) is proved.

$(3)$ If $u\not\equiv0$ is a critical point of $I$, then $I'(u)=0$
and thus $G(u)=0$. So $u$ is a critical point of $I$ constrained on
$\mathcal{N}$. Conversely, if $u$ is a critical point of $I$
constrained on $\mathcal{N}$, then there exists $\lambda\in\emph{R}$
such that $I'(u)=\lambda G'(u)$. It follows from  $u\in\mathcal{N}$
that
$$
\langle \lambda G'(u),u\rangle=\langle I'(u),u\rangle=0,
$$
which combined with \eqref{eq:<G'(u),u>} gives that $\lambda=0$.
Thus, $I'(u)=0$. So (3) is proved. $\quad\Box$

We recall some known facts. The positive solutions of the following problem
\begin{equation}\label{eq:Sch-limit}
  \left\{
\begin{array}{ll}
-\Delta u=u^{5}, \quad x\in\mathbb{R}^{3},\\
\\
u\in D^{1,2}(\mathbb{R}^{3})
\end{array}\right.
\end{equation}
must be of the form
\begin{equation}\label{eq:varphi_delta,x0}
\varphi_{\delta,y}=\frac{(3\delta)^{\frac{1}{4}}}{(\delta+|x-y|^{2})^{\frac{1}{2}}},\quad \delta>0,\quad y\in\mathbb{R}^{3}.
\end{equation}
The best Sobolev constant $S$ defined in (\ref{e1.1}) is achieved by the form \eqref{eq:varphi_delta,x0}, and
\begin{equation}\label{e3.2}
S=\big|\varphi_{\delta,y}\big|^{4}_{L^{6}(\mathbb{R}^{3})}=\big\|\varphi_{\delta,y}\big\|^{4/3}_{\mathbb{R}^{3}}.
\end{equation}
Moreover, we remark that the limit equation of \eqref{eq:Omega} is the following Kirchhoff type problem
\begin{equation}\label{eq:limit}
  \left\{
\begin{array}{ll}
\mathlarger{-\left(a+b\int_{\mathbb{R}^{3}}|\nabla u|^{2}dx\right)\Delta u=u^{5}}, \quad x\in\mathbb{R}^{3},\\
\\
u\in D^{1,2}(\mathbb{R}^{3}),\\
\end{array}\right.
\end{equation}
and the associated energy functional  $\Psi: D^{1,2}(\mathbb{R}^{3})\to \mathbb R$  is given by
\begin{equation}\label{eq:Psi(u)}
\Psi(u):=\frac{a}{2}\int_{\mathbb{R}^{3}}|\nabla u|^{2}dx+\frac{b}{4}\left(\int_{\mathbb{R}^{3}}|\nabla
u|^{2}dx\right)^{2}-\frac{1}{6}\int_{\mathbb{R}^{3}}u^{6}dx.
\end{equation}
Denote the least energy of the functional $\Psi$ by $M_{\infty}$, where
$$
M_{\infty}:=\inf\Big\{\Psi(u):\ u\in \mathcal{N}_{\infty}\Big\},\quad\mathcal{N}_{\infty}:=\{u\in D^{1,2}(\mathbb{R}^{3})\backslash\{0\}:\langle \Psi'(u),u\rangle=0\}.
$$
Then a ground state solution of \eqref{eq:limit} is a nontrivial solution $u\in D^{1,2}(\mathbb{R}^3)$ satisfying
$$\Psi(u)=M_{\infty} \ \ \text{and} \ \ \Psi'(u)=0.$$

In the next section we would analyze the behavior of a Palais-Smale sequence of $I(u)$ and give some compactness results which are
necessary in our study. Some properties about the ground state solution to the equation \eqref{eq:limit} are key in studying the global compactness in the usual sense, so we recall the following lemma.
\\
\textbf{Lemma 2.2.} (\cite{Xie-1}) {\it Let $a>0$, $b\geq 0$. Then the problem \eqref{eq:limit} has a positive ground state solution, unique up to translation and scaling.
Let $\varphi_{\delta,y}$ be defined in \eqref{eq:varphi_delta,x0}, then  the positive ground state solution  $\overline{u}_{\delta,y}$ of \eqref{eq:limit} is of  the form
\begin{equation}\label{overline{u}_delta,y}
\overline{u}_{\delta,y}(x):=\sqrt K\varphi_{\delta,y}(x)=\frac{\sqrt K(3\delta)^{\frac{1}{4}}}{(\delta+|x-y|^{2})^{\frac{1}{2}}},\quad \delta>0,\ y\in\mathbb{R}^{3},
\end{equation}
where  $K=\frac{1}{2}\left(bS^{\frac{3}{2}}+\sqrt{b^{2}S^{3}+4a}\right)$, $S$ is the best Sobolev constant.  And the least energy of the functional $\Psi$ is:
\begin{equation}\label{eq:M_infty}
M_{\infty}=\frac{a}{3}KS^{\frac{3}{2}}+\frac{b}{12}K^{2}S^{3}.
\end{equation}
Moreover, if $w$ is a sign-changing solution of the problem \eqref{eq:limit}, then $\Psi(w)>2M_{\infty}$.}\\
\textbf{Lemma 2.3.} (\cite{Willem}) {\it If $N\geq3$ and $h\in L^{\frac{N}{2}}(\mathbb{R}^{N})$, the functional $\chi:D^{1,2}(\mathbb{R}^{N})\rightarrow\mathbb{R}^{N}$ defined by
$$
\chi(u)=
\int_{\mathbb{R}^{N}}h(x)u^{2}dx
$$
is weakly continuous.}

We define
$$
M_{\Omega}:=\inf\Big\{I(u):\ u\in \mathcal{N}\Big\}.
$$
Then a ground state solution of \eqref{eq:Omega} is a nontrivial solution $u\in D_0^{1,2}(\Omega)$ which satisfies
$$
I(u)=M_{\Omega} \ \ \text{and} \ \  I'(u)=0.
$$
We have the following lemma.\\
\textbf{Lemma 2.4.} {\it Let $a>0$, $b\geq 0$, $V\in L^{3/2}(\Omega)$ be a non-negative function. Then $M_{\Omega}=M_{\infty}$ and $M_{\Omega}$ is not achieved; that is, the equation \eqref{eq:Omega} has no ground state solution.}\\
\textbf{Proof.} (1) We show that
\begin{equation}\label{eq:ground state}
M_{\Omega}=M_{\infty}.
\end{equation}
For $\forall u\in D_{0}^{1,2}(\Omega)$, let $u\equiv0$ outside
$\Omega$, then it can be extended to $D^{1,2}(\mathbb{R}^{3})$. Thus,
we may consider $D_{0}^{1,2}(\Omega)$ as a subspace of
$D^{1,2}(\mathbb{R}^{3})$. So for any $u\in\mathcal{N}$, there exists $\xi_{u}>0$ such that $\xi_{u}u\in\mathcal{N}_{\infty}$ and
$$\begin{array}{rcl}
M_{\infty}&\leq& \Psi(\xi_{u}u)=\mathlarger{\frac{a}{2}\|\xi_{u}u\|_{\mathbb{R}^{3}}^{2}+\frac{b}{4}\|\xi_{u}u\|_{\mathbb{R}^{3}}^{4}-\frac{1}{6}|\xi_{u}u|_{L^{6}(\mathbb{R}^{3})}^{6}}\\
\\
&\leq&\mathlarger{\frac{a}{2}\|\xi_{u}u\|_{\mathbb{R}^{3}}^{2}+\frac{1}{2}\int_{\Omega}V(x)(\xi_{u}u)^{2}dx+\frac{b}{4}\|\xi_{u}u\|_{\mathbb{R}^{3}}^{4}-\frac{1}{6}|\xi_{u}u|_{L^{6}(\mathbb{R}^{3})}^{6}}\\
\\
&=&I(\xi_{u}u)\leq I(u),
\end{array}
$$
which combined with the  arbitrariness of $u\in \mathcal{N}$ gives that
\begin{equation}\label{eq:Omega>R3}
M_{\Omega}\geq M_{\infty}.
\end{equation}

Now consider the sequence $\{\phi_{n}\}\subset D_{0}^{1,2}(\Omega)$
defined by $\phi_{n}(x):=\zeta(x)\overline{u}_{n}$, where $\overline{u}_{n}(\cdot)=\overline{u}(\cdot-y_{n})$ and $\overline{u}=\overline{u}_{1,0}\in
D^{1,2}(\mathbb{R}^{3})$ defined in (\ref{overline{u}_delta,y}) is a positive solution of \eqref{eq:limit}, $\{y_{n}\}\subset\Omega$ is a sequence of points such that
$|y_{n}|\rightarrow+\infty$ as $n\rightarrow+\infty$, $\zeta:\mathbb{R}^{3}\rightarrow[0,1]$ is
defined by
$$
\zeta(x)=\widetilde{\zeta}\left(\frac{|x|}{\overline{\rho}}\right),\quad \overline{\rho}:=\inf\left\{\rho:\mathbb{R}^{3}\setminus\Omega\subset\overline{B_{\rho}(0)}\right\},
$$
and $\widetilde{\zeta}(t):\emph{R}^{+}\cup\{0\}\rightarrow[0,1]$ is a
non-decreasing function such that $\widetilde{\zeta}=0$ for
$\forall t\leq1$, and $\widetilde{\zeta}=1$ for $\forall t\geq2$.
Nextly, we prove that
\begin{equation}\label{eq:I(phi_{n})}
I(\phi_{n})\rightarrow M_{\infty}\ \ \text{and}\ \ \langle I'(\phi_{n}),\phi_{n}\rangle\rightarrow0,\ \ \text{as}\ \ n\to+\infty.
\end{equation}
It follows from the definition of $\overline{u}_{n}$ that $\|\overline{u}_{n}\|_{\mathbb{R}^3}=\|\overline{u}\|_{\mathbb{R}^3}$ and $\overline{u}_{n}\rightharpoonup0$  in $D^{1,2}(\mathbb{R}^{3})$ as $n\rightarrow +\infty$. Then by Lemma 2.3 we have
$$
\int_{\Omega}V(x)\phi_{n}^{2}dx\rightarrow0\quad \text{as} \ \ n\rightarrow+\infty.
$$
Moreover, by \eqref{overline{u}_delta,y}, if $n$ is large enough, we can give the estimates
$$\begin{array}{rcl}
\left|\mathlarger{\int}_{\Omega}\phi_{n}^{6}dx-\mathlarger{\int}_{\mathbb{R}^{3}}\overline{u}^{6}dx\right|&
=&\mathlarger{\int}_{B_{2\bar{\rho}}}|\zeta(x)\overline{u}(x-y_{n})|^{6}dx
\\ \\
&\leq&k_{1}\mathlarger{\int_{B_{2\bar{\rho}}}\left(\frac{1}{|x-y_{n}|}\right)^{6}dx}=o\left(\frac{1}{|y_{n}|}\right)
\end{array}
$$
and
$$\begin{array}{rcl}
\big\|\phi_{n}(x)-\overline{u}(x-y_{n})\big\|^{2}_{\mathbb{R}^{3}}
&\leq&k_{2}\mathlarger{\int}_{B_{2\bar{\rho}}}\big|\nabla\overline{u}(x-y_{n})\big|^{2}dx
+k_{3}\mathlarger{\int}_{B_{2\bar{\rho}}}\big|\overline{u}(x-y_{n})\big|^{2}dx\\
\\
&=&o\left(\frac{1}{|y_{n}|}\right),
\end{array}
$$
where $k_{i}$ $(i=1,2,3)$ are positive constants. Thus, \eqref{eq:I(phi_{n})} is proved.

For $\phi_{n}\in D^{1,2}_{0}(\Omega)$, by the properties of
$\mathcal{N}$ in Lemma 2.1, we see that there exists a unique
$t_{n}\in(0,+\infty)$ such that $t_{n}\phi_{n}\in\mathcal{N}$; that
is,
$$
\langle I'(t_{n}\phi_{n}), t_{n}\phi_{n}\rangle=0.
$$
Now we show that
\begin{equation}\label{tn rightarrow1}
t_{n}\rightarrow1 \ \  \text{as}\  \  n\rightarrow+\infty.
\end{equation}
In fact, by the definition of $\phi_{n}$, there exist constants
$c_{1},c_{2}>0$ such that $c_{1}\leq\|\phi_{n}\|,|\phi_{n}|_6\leq c_{2}$, which
combined with \eqref{eq:u lower bound} gives that $t_{n}\geq t_{0}$
for some $t_{0}>0$. By  $t_{n}\phi_{n}\in\mathcal{N}$, we have
$$
0=\langle I'(t_{n}\phi_{n}), t_{n}\phi_{n}\rangle
=at_{n}^{2}\|\phi_{n}\|^{2}+t_{n}^{2}\int_{\Omega}V(x)\phi_{n}^{2}dx
+bt_{n}^{4}\|\phi_{n}\|^{4}-t_{n}^{6}|\phi_{n}|_{6}^{6}.
$$
Suppose $t_{n}\rightarrow+\infty$, then we have
$$
b\|\phi_{n}\|^{4}=t_{n}^{2}|\phi_{n}|_{6}^{6}+o_n(1)\rightarrow+\infty,\ \text{as}\ n\to+\infty,
$$
which is a contradiction. Thus, $\{t_{n}\}$ is bounded from above. Then
according to the properties of $\alpha_{u}(t)$ (see the proof of Lemma 2.1) and
$$
\alpha'_{\phi_{n}}(t_{n})=\frac{1}{t_{n}}\langle I'(t_{n}\phi_{n}),
t_{n}\phi_{n}\rangle,\quad \langle
I'(\phi_{n}),\phi_{n}\rangle\rightarrow0,
$$
we obtain \eqref{tn rightarrow1}. Thus, by
\eqref{eq:I(phi_{n})}, (\ref{tn rightarrow1}) and the fact $\int_{\Omega}V(x)\phi_{n}^{2}dx\rightarrow0$, we deduce that $I(t_{n}\phi_{n})\rightarrow
M_{\infty}$, which combined with
$t_{n}\phi_{n}\in\mathcal{N}$ gives that $M_{\Omega}\leq M_{\infty}$. This combined with  \eqref{eq:Omega>R3}
gives  \eqref{eq:ground state}.

(2) By contradiction, suppose that there exists a $v\in D_{0}^{1,2}(\Omega)$ such that
$$
I(v)=M_{\Omega}=M_{\infty}\ \text{and}\  I'(v)=0.
$$
By putting $v\equiv0$ in $\mathbb{R}^{3}\setminus\Omega$, $v$
could be regarded as an element of $D^{1,2}(\mathbb{R}^{3})$. Let $t_{0}>0$ be such that $t_{0}v\in\mathcal{N}_{\infty}$, then we have
$$\begin{array}{rcl}
M_{\infty}&\leq& \Psi(t_{0}v)=\mathlarger{\frac{a}{2}\|t_{0}v\|_{\mathbb{R}^{3}}^{2}+\frac{b}{4}\|t_{0}v\|_{\mathbb{R}^{3}}^{4}-\frac{1}{6}|t_{0}v|_{L^{6}(\mathbb{R}^{3})}^{6}}\\
\\
&\leq&\mathlarger{\frac{a}{2}\|t_{0}v\|_{\mathbb{R}^{3}}^{2}+\frac{1}{2}\int_{\Omega}V(x)(t_{0}v)^{2}dx+\frac{b}{4}\|t_{0}v\|_{\mathbb{R}^{3}}^{4}-\frac{1}{6}|t_{0}v|_{L^{6}(\mathbb{R}^{3})}^{6}}\\
\\
&=&I(t_{0}v)\leq I(v)=M_{\infty},
\end{array}
$$
which implies
$$
t_{0}=1\ \ \text{and}\ \  \int_{\Omega}V(x)v^{2}dx=0;
$$
that is, $v$ would be a minimizer of $M_{\infty}$. By the maximum principle, $v$ is strictly positive in $\mathbb{R}^{3}$,
which is a contradiction. The proof is complete. $\quad\Box$
\\
\textbf{Lemma 2.5.}(\cite{Willem}) {\it Let $N\geq3$. If $u_{n}\rightharpoonup u$ in $D^{1,2}(\mathbb{R}^{N})$ and $u\in L_{loc}^{\infty}(\mathbb{R}^{N})$, then}
$$
|u_{n}|^{2^{*}-2}u_{n}-|u_{n}-u|^{2^{*}-2}(u_{n}-u)\rightarrow |u|^{2^{*}-2}u \quad \text{in} \ (D^{1,2}(\mathbb{R}^{N}))'.
$$
\textbf{Lemma 2.6.}(\cite{{Benci-1}}) {\it Let $N\geq 3$, $\{u_{n}\}\subset H^{1}_{loc}(\mathbb{R}^{N})$ be a sequence of functions such that
$$
u_{n}\rightharpoonup0\quad\text{in} \ \ H^{1}(\mathbb{R}^{N}).
$$
Assume that there exist a bounded open set $Q\subset\mathbb{R}^{N}$ and a positive constant $\gamma$ such that
$$
\int_{Q}|\nabla u_{n}|^{2}\geq\gamma>0,\quad \int_{Q}|u_{n}|^{2^{*}}\geq\gamma>0.
$$
Moreover, suppose that
$$
\Delta u_{n}+|u_{n}|^{2^{*}-2}u_{n}=\chi_{n},
$$
where $\chi_{n}\in H^{-1}(\mathbb{R}^{N})$ and $\langle\chi_{n},\phi\rangle\leq\varepsilon_{n}\|\phi\|_{H^{1}}$, $\phi\in C_{0}^{\infty}(U)$, $U$ is an open neighborhood of $Q$ and $\varepsilon_{n}$ converges to 0 as $n\rightarrow+\infty$. Then there exist a sequence of positive numbers $\{\sigma_{n}\}$ and a sequence of points $\{y_{n}\}\subset\mathbb{R}^{N}$ such that
$$
v_{n}(x):=\sigma_{n}^{\frac{N-2}{2}}u_{n}(\sigma_{n}x+y_{n})
$$
converges weakly in $D^{1,2}(\mathbb{R}^{N})$ to a nontrivial solution $v$ of $-\Delta u=|u|^{2^{*}-2}u$. Moreover,}
$$
y_{n}\rightarrow y\in\overline{Q}\ \text{and}\  \sigma_{n}\rightarrow0.
$$

\section{Compactness lemma}

In order to have a better understanding of how the Palais-Smale condition may fail, we need to investigate more closely the compactness problem in this section. We give a slightly more general result, which should be useful in other contexts. To this end, we define a functional with critical growth term by
\begin{equation}\label{eq:I(u)-R^N}
I_{*}(u):=\frac{a}{2}\int_{\Omega}|\nabla u|^{2}dx+\frac{b}{4}\left(\int_{\Omega}|\nabla
u|^{2}dx\right)^{2}+\frac{1}{2}\int_{\Omega}V(x)u^{2}dx-\frac{1}{2^{*}}\int_{\Omega}|u|^{2^{*}}dx,
\end{equation}
where $N\geq 3$, $\Omega\subset\mathbb{R}^{N}$, $2^{*}:=\frac{2N}{N-2}$.
\\
\textbf{Lemma 3.1.} {\it Assume that $N\geq 3$, $a>0$, $b\geq 0$ and $V\in L^{N/2}(\Omega)$. Let $\{u_{n}\}\subset D_{0}^{1,2}(\Omega)$ be a sequence such that
\begin{equation}\label{eq:PS-sequence}
I_{*}(u_{n})\rightarrow c\  \ \text{and}\  \ I_{*}'(u_{n})\rightarrow0,\  \ \text{as}\ \ n\to+\infty.
\end{equation}
Then there exist a number $k\in\mathbb{N}:=\{0,1,2,\cdots\}$, $k$ sequences of numbers $\{\sigma_{n}^{i}\}_{n\in\mathbb{N}}\subset\mathbb{R}^{+}$, points $\{y_{n}^{i}\}_{n\in\mathbb{N}}\subset\mathbb{R}^N$, $1\leq i\leq k$, $k+1$ sequences of functions $\{u_n^{(j)}\}_n\subset D^{1,2}(\mathbb{R}^{N})$, $0\leq j\leq k$, such that for some subsequence, still denoted by $\{u_{n}\}$,
\begin{equation}\label{eq:u(n)-decomposition}
\begin{array}{rcl}
&[i]&\quad u_{n}(x)=u^{(0)}_{n}(x)+\mathlarger{\sum_{i=1}^{k}(\sigma_{n}^{i})^{-\frac{N-2}{2}}u_n^{(i)}\left(\frac{\cdot-y_{n}^{i}}{\sigma_{n}^{i}}\right)},\\
\\
&[ii]&\quad u^{(0)}_{n}(x)\rightarrow u^{(0)}(x)\  \text{in} \ D_{0}^{1,2}(\Omega) \ \text{as}\   n\rightarrow+\infty,\\
\\
&[iii]&\quad u^{(i)}_{n}(x)\rightarrow u^{(i)}(x) \ \text{in} \ D^{1,2}(\mathbb{R}^{N})\ \text{as}\   n\rightarrow+\infty,\ 1\leq i\leq k,\\
\end{array}
\end{equation}
and $u^{(0)}$, $u^{(i)}$ ($1\leq i\leq k$) satisfy
\begin{equation}\label{eq:u(0)-function}
-\left(a+bA^{2}\right)\Delta u^{(0)}+V(x)u^{(0)}=|u^{(0)}|^{2^{*}-2}u^{(0)}, \quad in \ \Omega,
\end{equation}
\begin{equation}\label{eq:u(i)-function}
-\left(a+bA^{2}\right)\Delta u^{(i)}=|u^{(i)}|^{2^{*}-2}u^{(i)}, \quad in \ \mathbb{R}^{N},\ 1\leq i\leq k,
\end{equation}
where $A\in [0,+\infty)$ satisfies
\begin{equation}\label{e2.2}
A^{2}=|\nabla u^{(0)}|^{2}_{2}+\sum_{i=1}^{k}|\nabla u^{(i)}|^{2}_{L^{2}(\mathbb{R}^{N})}.
\end{equation}
Moreover, as $n\rightarrow+\infty$, we have
\begin{equation}\label{eq:un-decomposition}
\left\|u_{n}-u^{(0)}-\sum_{i=1}^{k}(\sigma_{n}^{i})^{-\frac{N-2}{2}}u^{(i)}\left(\frac{\cdot-y_{n}^{i}}{\sigma_{n}^{i}}\right)\right\|\rightarrow0,\quad \|u_{n}\|^{2}\rightarrow A^{2},
\end{equation}
and
\begin{equation}\label{eq:I(u)-decomposition}
I_{*}(u_{n})\rightarrow
\widetilde{I}(u^{(0)})+\sum_{i=1}^{k}\widetilde{\Psi}(u^{(i)}),
\end{equation}
where
$$\begin{array}{rcl}
\widetilde{I}(u)&:=&\mathlarger{\left(\frac{a}{2}+\frac{bA^{2}}{4}\right)\|u\|^{2}+\frac{1}{2}\int_{\Omega}V(x)u^{2}dx-\frac{1}{2^{*}}\int_{\Omega}|u|^{2^{*}}dx},\\
\\
\widetilde{\Psi}(u)&:=&\mathlarger{\left(\frac{a}{2}+\frac{bA^{2}}{4}\right)\|u\|_{\mathbb{R}^{N}}^{2}-\frac{1}{2^{*}}\int_{\mathbb{R}^{N}}|u|^{2^{*}}dx}.
\end{array}
$$}

For any non-negative constant $A$, we define two functionals as follows:
$$\begin{array}{rcl}
\widetilde{I}_{*}(u)&:=&\mathlarger{\frac{a+bA^{2}}{2}\int_{\Omega}|\nabla u|^{2}dx+\frac{1}{2}\int_{\Omega}V(x)u^{2}dx-\frac{1}{2^{*}}\int_{\Omega}|u|^{2^{*}}dx},\\
\\
\widetilde{\Psi}_{*}(u)&:=&\mathlarger{\frac{a+bA^{2}}{2}\int_{\mathbb{R}^{N}}|\nabla u|^{2}dx-\frac{1}{2^{*}}\int_{\mathbb{R}^{N}}|u|^{2^{*}}dx}.
\end{array}
$$
To prove Lemma 3.1, we first give the following lemma. \\
\textbf{Lemma 3.2.} {\it Let $\{u_{n}\}\subset D^{1,2}(\mathbb{R}^N)$ be a Palais-Smale sequence of $\widetilde{\Psi}_{*}$ such that $u_{n}\in C_{0}^{\infty}(\mathbb{R}^{N})$,
$$
u_{n}\rightharpoonup0\ \ \text{in} \ \ D^{1,2}(\mathbb{R}^{N})\ \  \text{and}\ \  u_{n}\not\rightarrow0\ \ \text{in} \ \  D^{1,2}(\mathbb{R}^{N}),\  \ \text{as}\ \  n\to+\infty.
$$
Then there exist a sequence of positive numbers $\{\sigma_{n}\}$ and a sequence of points $\{y_{n}\}\subset\mathbb{R}^{N}$ such that
$$
v_{n}(x):=\sigma_{n}^{\frac{N-2}{2}}u_{n}(\sigma_{n}x+y_{n})
$$
converges weakly in $D^{1,2}(\mathbb{R}^{N})$ to a nontrivial solution $v$ of the problem \eqref{eq:u(i)-function}.}
\\
\textbf{Proof.} By using Lemma 2.6 and arguing as in the proof of \cite[Lemma 5.3]{Xie-1},  we can prove this lemma and the details will be omitted here. $\quad\Box$

Now we give the proof of Lemma 3.1.\\
\textbf{Proof of Lemma 3.1.}  First we prove that  $\{u_{n}\}$ is bounded in $D_{0}^{1,2}(\Omega)$. By \eqref{eq:PS-sequence} we have
$$\begin{array}{cl}
c+o_n(1)+o_n(1)\|u_n\|&=I_*(u_{n})-\mathlarger{\frac{1}{2^{*}}}\big\langle I'_*(u_{n}),u_{n}\big\rangle\\
\\
&=\mathlarger{(\frac{1}{2}-\frac{1}{2^{*}})\|u_{n}\|^{2}+(\frac{1}{2}-\frac{1}{2^{*}})\int_{\Omega}V(x)u_{n}^{2}dx}
+\mathlarger{(\frac{1}{4}-\frac{1}{2^{*}})\|u_{n}\|^{4}},\\
\end{array}
$$
then by $2^{*}=\frac{2N}{N-2}$ and $N\geq3$, it is easy to see $\{u_{n}\}$ is bounded in $D_{0}^{1,2}(\Omega)$ when $N=3,4$. When $N\geq5$, we assume $\|u_{n}\|\rightarrow+\infty$ by contradiction, then
$$
0<o_n(1)+(\frac{1}{2^{*}}-\frac{1}{4})\|u_{n}\|^{2}\leq C_{1}<\infty,
$$
which is a contradiction. Thus $\{u_{n}\}$ is bounded in $D_{0}^{1,2}(\Omega)$ for $N\geq3$.
So there exist $u^{(0)}\in D_{0}^{1,2}(\Omega)$ and $A\in\emph{R}$ such that up to a subsequence,
still denoted by $\{u_n\}$,
\begin{equation}\label{eq:u_{n}-weakly}
\begin{array}{cl}
&u_{n}\rightharpoonup u^{(0)} \  \ \text{in} \  \ D_{0}^{1,2}(\Omega),\\
\\
&u_{n}(x)\rightarrow u^{(0)}(x)\ \  \text{a.e.\ in} \ \ \Omega,\\
\end{array}
\end{equation}
and
\begin{equation}\label{eq:A^2}
\int_{\Omega}|\nabla u_{n}|^{2}dx\rightarrow A^{2}
\end{equation}
as $n\to +\infty$.
By standard arguments, we have $\widetilde{I}'_{*}(u^{(0)})=0$. Consequently, $u^{(0)}$ solves \eqref{eq:u(0)-function}. Moreover, by using an iteration argument, we could deduce that
\begin{equation}\label{eq:u(0)-L-loc}
u^{(0)}\in L^{\infty}_{loc}(\mathbb{R}^{N}).
\end{equation}
Now let
\begin{equation}\label{eq:v_{n}^{(1)}}
v_{n}^{(1)}(x):=\left\{
\begin{array}{ll}
(u_{n}-u^{(0)})(x), \quad x\in\Omega,\\
\\
0,\ \qquad\qquad\qquad x\in\mathbb{R}^{N}\setminus\Omega.\\
\end{array}\right.
\end{equation}
By \eqref{eq:u_{n}-weakly}, $v_{n}^{(1)}\rightharpoonup0$ in $D^{1,2}(\mathbb{R}^{N})$ as $n\rightarrow+\infty$, then it follows from Lemma 2.3 that
\begin{equation}\label{eq:V(x)vn(1)}
\int_{\mathbb{R}^{N}}V(x)(v_{n}^{(1)})^2dx\rightarrow0\ \ \text{as}\ \ n\rightarrow+\infty.
\end{equation}
Thus,
\begin{equation}\label{eq:vn(1)-decomposition}
\begin{array}{cl}
&\|v_{n}^{(1)}\|^{2}_{\mathbb{R}^{N}}=\|v_{n}^{(1)}\|^{2}=\|u_{n}\|^{2}-\|u^{(0)}\|^{2}+o_n(1),\\
\\
&\widetilde{\Psi}_{*}(v_{n}^{(1)})=\widetilde{I}_{*}(v_{n}^{(1)})+o_n(1)=\widetilde{I}_{*}(u_{n})-\widetilde{I}_{*}(u^{(0)})+o_n(1),\\
\\
&\widetilde{I}_{*}(u_{n})=\mathlarger{{I}_{*}(u_{n})+\frac{b}{4}A^4+o_n(1)=c+\frac{b}{4}A^4+o_n(1)}.
\end{array}
\end{equation}
Nextly, we will divide the remaining proof into several steps.\\
\textbf{Step 1}: We have two possibilities.\\
(1) If $v_{n}^{(1)}\rightarrow0$  in $D^{1,2}(\mathbb{R}^{N})$, the lemma is proved with $k=0$.\\
(2) If $v_{n}^{(1)}\not\rightarrow0$  in $D^{1,2}(\mathbb{R}^{N})$, by \eqref{eq:u(0)-L-loc} and Lemma 2.5, for any $\varphi\in C_{0}^{\infty}(\mathbb{R}^{N})$, we have
$$
\int_{\mathbb{R}^{N}}\left(|v_{n}^{(1)}|^{2^{*}-2}v_{n}^{(1)}+|u^{(0)}|^{2^{*}-2}u^{(0)}-|u_{n}|^{2^{*}-2}u_{n}\right)\varphi dx=o_{n}(1)\|\varphi\|_{\mathbb{R}^N}.
$$
Then it follows from \eqref{eq:PS-sequence} and $\widetilde{I}'_{*}(u^{(0)})=0$ that
\begin{equation}\label{eq:widetilde-vn(1)-decomposition}
\begin{array}{cl}
\langle\widetilde{\Psi}'_{*}(v_{n}^{(1)}),\varphi\rangle&=\langle\widetilde{I}'_{*}(v_{n}^{(1)}),\varphi\rangle+o_{n}(1)\|\varphi\|_{\mathbb{R}^N}\\
\\
&=o_{n}(1)\|\varphi\|_{\mathbb{R}^N}+\langle\widetilde{I}'_{*}(u_{n}),\varphi\rangle-\langle\widetilde{I}'_{*}(u^{(0)}),\varphi\rangle\\
\\
&\qquad-\mathlarger{\int}_{\mathbb{R}^{N}}\left(|v_{n}^{(1)}|^{2^{*}-2}v_{n}^{(1)}+|u^{(0)}|^{2^{*}-2}u^{(0)}-|u_{n}|^{2^{*}-2}u_{n}\right)\varphi dx\\
\\
&=\langle I_{*}'(u_{n}),\varphi\rangle-\langle\widetilde{I}'_{*}(u^{(0)}),\varphi\rangle+o_{n}(1)\|\varphi\|_{\mathbb{R}^N}
=o_{n}(1)\|\varphi\|_{\mathbb{R}^N}.
\end{array}
\end{equation}
Thus, $\widetilde{\Psi}'_{*}(v_{n}^{(1)})\rightarrow0$ as $n\rightarrow+\infty$; that is, $\{v_{n}^{(1)}\}$ is a Palais-Smale sequence for $\widetilde{\Psi}_{*}$. By the definition of $D^{1,2}(\mathbb{R}^{N})$, for any $n$, there exists $u_{n}^{(1)}\in C_{0}^{\infty}(\mathbb{R}^{N})$ such that
$$
\big\|u_{n}^{(1)}-v_{n}^{(1)}\big\|_{\mathbb{R}^{N}}<\frac{1}{n}\ \  \text{and}\  \ \big|\widetilde{\Psi}'_{*}(u_{n}^{(1)})-\widetilde{\Psi}'_{*}(v_{n}^{(1)})\big|<\frac{1}{n}.
$$
Then it is easy to see that
\begin{equation}\label{eq:u_n(1)}
\begin{array}{cl}
&\big\|u_{n}^{(1)}\big\|_{\mathbb{R}^{N}}^{2}=\big\|v_{n}^{(1)}\big\|_{\mathbb{R}^{N}}^{2}+o_{n}(1)=\|u_{n}\|^{2}-\|u^{(0)}\|^{2}+o_n(1),\\
\\
&\widetilde{\Psi}_{*}(u_{n}^{(1)})=\widetilde{\Psi}_{*}(v_{n}^{(1)})+o_{n}(1)=\widetilde{I}_{*}(u_{n})-\widetilde{I}_{*}(u^{(0)})+o_n(1),\\
\\
&\widetilde{\Psi}'_{*}(u_{n}^{(1)})=\widetilde{\Psi}'_{*}(v_{n}^{(1)})+o_{n}(1)=o_{n}(1).
\end{array}
\end{equation}
Thus, $\{u_{n}^{(1)}\}$ is a Palais-Smale sequence for $\widetilde{\Psi}_{*}$ and $\{u_{n}^{(1)}\}\subset C_{0}^{\infty}(\mathbb{R}^{N})$,
\begin{equation}\label{eq:u{n}(1)-weakly}
\begin{array}{cl}
&u_{n}^{(1)}\rightharpoonup0  \ \text{in} \ D^{1,2}(\mathbb{R}^{N}),\\
\\
&u_{n}^{(1)}\not\rightarrow0 \ \text{in} \ D^{1,2}(\mathbb{R}^{N}).\\
\end{array}
\end{equation}
Then by Lemma 3.2, there exist a sequence of points $\{y_{n}^{1}\}\subset\mathbb{R}^{N}$ and a sequence of positive numbers $\{\sigma_{n}^{1}\}$ such that
\begin{equation*}
\begin{array}{cl}
&\widetilde{u}_{n}^{(1)}(\cdot):=(\sigma_{n}^{1})^{\frac{N-2}{2}}u_{n}^{(1)}(\sigma_{n}^{1}\cdot+y_{n}^{1}),\\
\\
&\widetilde{u}_{n}^{(1)}\rightharpoonup u^{(1)} \   \text{in} \ D^{1,2}(\mathbb{R}^{N}),\\
\\
&\widetilde{\Psi}'_{*}(u^{(1)})=0,\quad u^{(1)}\neq0.\\
\end{array}
\end{equation*}
Thus, $u^{(1)}$ is a nontrivial solution of the problem \eqref{eq:u(i)-function}. And
$$
\widetilde{\Psi}_{*}(u_{n}^{(1)})=\widetilde{\Psi}_{*}(\widetilde{u}_{n}^{(1)})=\widetilde{\Psi}_{*}(u^{(1)})+\widetilde{\Psi}_{*}(\widetilde{u}_{n}^{(1)}-u^{(1)})+o_{n}(1),
$$
$$
\big\|u_{n}^{(1)}\big\|_{\mathbb{R}^{N}}^{2}=\big\|\widetilde{u}_{n}^{(1)}\big\|_{\mathbb{R}^{N}}^{2}
=\big\|u^{(1)}\big\|_{\mathbb{R}^{N}}^{2}+\big\|\widetilde{u}_{n}^{(1)}-u^{(1)}\big\|_{\mathbb{R}^{N}}^{2}+o_{n}(1),
$$
which combined with  \eqref{eq:u_n(1)} gives  that
$$
\widetilde{I}_{*}(u_{n})=\widetilde{I}_{*}(u^{(0)})+\widetilde{\Psi}_{*}(u^{(1)})+\widetilde{\Psi}_{*}(\widetilde{u}_{n}^{(1)}-u^{(1)})+o_{n}(1),
$$
$$
\big\|u_{n}\big\|^{2}=\big\|u^{(0)}\big\|^{2}+\big\|u^{(1)}\big\|_{\mathbb{R}^{N}}^{2}+\big\|\widetilde{u}_{n}^{(1)}-u^{(1)}\big\|_{\mathbb{R}^{N}}^{2}+o_{n}(1).
$$
\textbf{Step 2}: Let $v_{n}^{(2)}:=\widetilde{u}_{n}^{(1)}-u^{(1)}$. We have two possibilities.\\
(1) If $v_{n}^{(2)}\rightarrow 0$  in $D^{1,2}(\mathbb{R}^{N})$, the lemma is proved with $k=1$.\\
(2) If $v_{n}^{(2)}\not\rightarrow 0$  in $D^{1,2}(\mathbb{R}^{N})$, then, analogously $v_{n}^{(2)}$ is a Palais-Smale sequence of $\widetilde{\Psi}_{*}$.
So we can find $\{u_{n}^{(2)}\}\subset C_{0}^{\infty}(\mathbb{R}^{N})$ such that
$$
\big\|u_{n}^{(2)}-v_{n}^{(2)}\big\|_{\mathbb{R}^{N}}<\frac{1}{n}\ \text{and}\  \big|\widetilde{\Psi}'_{*}(u_{n}^{(2)})-\widetilde{\Psi}'_{*}(v_{n}^{(2)})\big|<\frac{1}{n}.
$$
Then
\begin{equation}\label{e1.2}
\begin{split}
\|u_n^{(2)}\|_{\mathbb{R}^N}^2=\|v_n^{(2)}\|_{\mathbb{R}^N}^2+o_n(1)&=\|\widetilde{u}_n^{(1)}-u^{(1)}\|_{\mathbb{R}^N}^2+o_n(1)\\
&=\|u_n\|^2-\|u^{(0)}\|^2-\|u^{(1)}\|_{\mathbb{R}^N}^2+o_n(1),\\
\widetilde{\Psi}_{*}(u_{n}^{(2)})=\widetilde{\Psi}_{*}(v_{n}^{(2)})+o_n(1)&=\widetilde{\Psi}_{*}(\widetilde{u}_{n}^{(1)}-u^{(1)})+o_n(1)\\
&=\widetilde{I}_{*}(u_n)-\widetilde{I}_{*}(u^{(0)})-\widetilde{\Psi}_{*}(u^{(1)})+o_n(1),\\
\widetilde{\Psi}'_{*}(u_{n}^{(2)})&=\widetilde{\Psi}_{*}'(v_{n}^{(2)})+o_n(1)=o_n(1).
\end{split}
\end{equation}
Thus, $\{u_n^{(2)}\}$ is a Palais-Smale sequence for $\widetilde{\Psi}_{*}$ and (\ref{eq:u{n}(1)-weakly}) holds for $u_n^{(2)}$. So there exist a sequence of points $\{y_{n}^{2}\}\subset\mathbb{R}^{N}$ and a sequence of numbers $\{\sigma_{n}^{2}\}\subset\mathbb{R}^{+}$ such that
\begin{equation*}
\begin{array}{cl}
&\widetilde{u}_{n}^{(2)}(\cdot):=(\sigma_{n}^{2})^{\frac{N-2}{2}}u_{n}^{(2)}(\sigma_{n}^{2}\cdot+y_{n}^{2}),\\
\\
&\widetilde{u}_{n}^{(2)}\rightharpoonup u^{(2)} \   \text{in} \ D^{1,2}(\mathbb{R}^{N}),\\
\\
&\widetilde{\Psi}'_{*}(u^{(2)})=0,\quad u^{(2)}\neq0.\\
\end{array}
\end{equation*}
Thus, $u^{(2)}$ is a nontrivial solution of the problem \eqref{eq:u(i)-function}. Moreover, there hold
$$
\widetilde{\Psi}_{*}(u_{n}^{(2)})=\widetilde{\Psi}_{*}(u^{(2)})+\widetilde{\Psi}_{*}(\widetilde{u}_{n}^{(2)}-u^{(2)})+o_{n}(1),
$$
$$
\big\|u_{n}^{(2)}\big\|_{\mathbb{R}^{N}}^{2}=\big\|u^{(2)}\big\|_{\mathbb{R}^{N}}^{2}+\big\|\widetilde{u}_{n}^{(2)}-u^{(2)}\big\|_{\mathbb{R}^{N}}^{2}+o_{n}(1),
$$
which combined with (\ref{e1.2}) gives that
$$
\widetilde{I}_{*}(u_{n})=\widetilde{I}_{*}(u^{(0)})+\widetilde{\Psi}_{*}(u^{(1)})+\widetilde{\Psi}_{*}(u^{(2)})+\widetilde{\Psi}_{*}(\widetilde{u}_{n}^{(2)}-u^{(2)})+o_{n}(1),
$$
$$
\big\|u_{n}\big\|^{2}=\big\|u^{(0)}\big\|^{2}+\big\|u^{(1)}\big\|_{\mathbb{R}^{N}}^{2}+\big\|u^{(2)}\big\|_{\mathbb{R}^{N}}^{2}+\big\|\widetilde{u}_{n}^{(2)}-u^{(2)}\big\|_{\mathbb{R}^{N}}^{2}+o_{n}(1).
$$
Iterating the above procedure, we obtain sequences $\{\widetilde{u}_{n}^{(m-1)}\}$ in this way.\\
\textbf{Step m}: Let $v_{n}^{(m)}:=\widetilde{u}_{n}^{(m-1)}-u^{(m-1)}$. We also have two possibilities.\\
(1) If $v_{n}^{(m)}\rightarrow 0$  in $D^{1,2}(\mathbb{R}^{N})$, the lemma is proved with $k=m-1$.\\
(2) If $v_{n}^{(m)}\not\rightarrow 0$  in $D^{1,2}(\mathbb{R}^{N})$, arguing as before, we obtain:
\begin{enumerate}
  \item[(a)] sequences $\{u_{n}^{(i)}\}\subset C_{0}^{\infty}(\mathbb{R}^{N})$ for $i=1,2,...,m$ such that
$$
\big\|u_{n}^{(i)}-v_{n}^{(i)}\big\|_{\mathbb{R}^{N}}<\frac{1}{n}\ \text{and}\ \big|\widetilde{\Psi}'_{*}(u_{n}^{(i)})-\widetilde{\Psi}'_{*}(v_{n}^{(i)})\big|<\frac{1}{n}.
$$
  \item[(b)] sequences of points $\{y_{n}^{i}\}\subset\mathbb{R}^{N}$ and sequences of positive numbers $\{\sigma_{n}^{i}\}\subset\mathbb{R}^{+}$ such that
\begin{equation*}
\begin{array}{cl}
&\widetilde{u}_{n}^{(i)}:=(\sigma_{n}^{i})^{\frac{N-2}{2}}u_{n}^{(i)}(\sigma_{n}^{(i)}\cdot+y_{n}^{i}),\\
\\
&\widetilde{u}_{n}^{(i)}\rightharpoonup u^{(i)}  \ \text{in} \ D^{1,2}(\mathbb{R}^{N}),\\
\\
&\widetilde{\Psi}'_{*}(u^{(i)})=0,\quad u^{(i)}\neq0,\\
\end{array}
\end{equation*}
for all $i=1,2,...,m$, and the procedure continues.
\end{enumerate}
But in fact, at some step $k+1$, the first case must occur to stop the process and the lemma would be proved. To show this, we notice
that by induction from the above  procedure, for any $l\in\emph{N}$,
\begin{equation}\label{eq:vn(l)-decomposition}
\|v_{n}^{(l)}\|^{2}_{\mathbb{R}^{N}}=\|u_{n}\|^{2}-\|u^{(0)}\|^{2}-\sum_{i}^{l-1}\|u^{(i)}\|^{2}_{\mathbb{R}^{N}}+o_{n}(1).
\end{equation}
In view of that $u^{(i)}$ (for any $i\geq 1$) is a nontrivial critical point of $\widetilde{\Psi}_{*}$, we have
$$
\|u^{(i)}\|_{\mathbb{R}^N}\geq C,\ \text{where}\  C>0\ \text{is\ a\ positive\ constant},
$$
which combined with \eqref{eq:vn(l)-decomposition} gives that the process has to stop at some index $k\geq0$. Thus,
\begin{equation}\label{eq:widetilde I*(un)-decomposition}
\widetilde{I}_{*}(u_{n})=\widetilde{I}_{*}(u^{(0)})+\sum_{i=1}^{k}\widetilde{\Psi}_{*}(u^{(i)})+o_{n}(1),
\end{equation}
\begin{equation}\label{eq:un-decomposition-2}
\big\|u_{n}\big\|^{2}=\big\|u^{(0)}\big\|^{2}+\sum_{i=1}^{k}\big\|u^{(i)}\big\|_{\mathbb{R}^{N}}^{2}+o_{n}(1).
\end{equation}
Then by \eqref{eq:A^2}, \eqref{eq:widetilde I*(un)-decomposition} and \eqref{eq:un-decomposition-2}, we have
\begin{equation*}
\begin{array}{rcl}
I_{*}(u_{n})&=&\widetilde{I}(u_{n})+o_{n}(1)\\
\\
&=&\widetilde{I}_{*}(u_{n})-\mathlarger{\frac{bA^{2}}{4}}\|u_{n}\|^{2}+o_{n}(1)\\
\\
&=&\mathlarger{\widetilde{I}_{*}(u^{(0)})-\frac{bA^{2}}{4}\|u^{(0)}\|^{2}+\sum_{i=1}^{k}\left(\widetilde{\Psi}_{*}(u^{(i)})-\frac{bA^{2}}{4}\|u^{(i)}\|_{\mathbb{R}^{N}}^{2}\right)+o_{n}(1)}\\
\\
&=&\widetilde{I}(u^{(0)})+\mathlarger{\sum_{i=1}^{k}}\widetilde{\Psi}(u^{(i)})+o_{n}(1),
\end{array}
\end{equation*}
which proves \eqref{eq:I(u)-decomposition}. By direct calculations, we have
\begin{equation*}
\begin{split}
u_n&=\left(u^{(0)}+o_n(1)\right)+(\sigma_n^1)^{-\frac{N-2}{2}}\left(u^{(1)}+o_n(1)\right)\left(\frac{x-y_n^1}{\sigma_n^1}\right)\\
&\qquad+(\sigma_n^1\sigma_n^2)^{-\frac{N-2}{2}}\left(u^{(2)}+o_n(1)\right)\left(\frac{x-y_n^1-\sigma_n^1y_n^2}{\sigma_n^1\sigma_n^2}\right)\\
&\qquad+(\sigma_n^1\sigma_n^2\sigma_n^3)^{-\frac{N-2}{2}}\left(u^{(3)}+o_n(1)\right)
\left(\frac{x-y_n^1-\sigma_n^1y_n^2-\sigma_n^1\sigma_n^2y_n^3}{\sigma_n^1\sigma_n^2\sigma_n^3}\right)\\
&\qquad+\cdots\\
&\qquad+(\sigma_n^1\sigma_n^2\cdots\sigma_n^k)^{-\frac{N-2}{2}}\left(u^{(k)}+o_n(1)\right)
\left(\frac{x-y_n^1-\sigma_n^1y_n^2-\cdots-\sigma_n^1\sigma_n^2\cdots \sigma_n^{k-1} y_n^k}{\sigma_n^1\sigma_n^2\cdots\sigma_n^k}\right).
\end{split}
\end{equation*}
Rewrite the notations, we obtain  (\ref{eq:u(n)-decomposition})-(\ref{eq:u(i)-function}). And \eqref{eq:un-decomposition} follows from (\ref{eq:u(n)-decomposition}), \eqref{eq:un-decomposition-2} and
$$
\int_{\mathbb{R}^{N}}|\nabla v^{(i)}|^{2}dx=\int_{\mathbb{R}^{N}}|\nabla u^{(i)}|^{2}dx,\quad \text{where} \  \  v^{(i)}:=(\sigma_{n}^{i})^{-\frac{N-2}{2}}u^{(i)}\left(\frac{\cdot-y_{n}^{i}}{\sigma_{n}^{i}}\right).
$$
The proof is complete. $\quad\Box$

In what follows, we assume $N=3$. By Lemma 3.1, we  obtain the following two lemmas.\\
\textbf{Lemma 3.3.} {\it  Let $a>0$, $b\geq 0$, $\Omega$ be as in (\ref{eq:Omega}), $V\in L^{3/2}(\Omega)$ be non-negative and continuous. If $\{u_{n}\}\subset D_{0}^{1,2}(\Omega)$ is a positive function sequence such that
\begin{equation}\label{eq:Lemma 3.3}
I(u_{n})\rightarrow M_{\Omega} \ \ \text{and}\  \ \langle I'(u_{n}),u_{n}\rangle=0, \ \ \text{as}\ \  n\to+\infty,
\end{equation}
then
\begin{equation}\label{e3.1}
u_{n}(x)=w_{n}(x)+\bar{u}_{\delta_{n},y_{n}}(x),
\end{equation}
where $\{w_{n}\}\subset D^{1,2}(\mathbb{R}^{3})$ is a sequence converging strongly to $0$ in $D^{1,2}(\mathbb{R}^{3})$, and $\bar{u}_{\delta_{n},y_{n}}$ defined in (\ref{overline{u}_delta,y}) is the positive function realizing $M_{\infty}$.}\\
\textbf{Proof.} By \eqref{eq:Lemma 3.3}, $\{u_{n}\}\subset D_{0}^{1,2}(\Omega)$ is a minimizing sequence for $I|_{\mathcal{N}}$, then by the Ekeland variational principle (\cite{Willem}), there exists a sequence $\{v_{n}\}\subset D_{0}^{1,2}(\Omega)$, $v_{n}\in\mathcal{N}$ such that
\begin{equation}\label{eq:vn-minimizing sequence of N}
I(v_{n})\rightarrow M_{\Omega},\quad I'(v_{n})-\lambda_{n}G'(v_{n})\rightarrow0,\quad\|u_{n}-v_{n}\|\rightarrow0,
\end{equation}
as $n\rightarrow+\infty$, where $\lambda_{n}\in\emph{R}$. We may assume $v_{n}\geq0$.

Now we claim that $I'(v_{n})\rightarrow0$ as $n\rightarrow+\infty$.
In fact, taking the scalar product with $v_{n}$, we obtain
$$
o_n(1)\|v_n\|=\langle I'(v_{n}),v_{n}\rangle-\lambda_{n}\langle
G'(v_{n}),v_{n}\rangle.
$$
For $v_{n}\in\mathcal{N}$, by \eqref{eq:<G'(u),u>} and \eqref{eq:u
lower bound}, we have $\langle
G'(v_{n}),v_{n}\rangle<\widetilde{c}<0$. By
$I(v_{n})\rightarrow M_{\Omega}$, $v_{n}\in\mathcal{N}$, the boundedness of $\{v_{n}\}$ could be
obtained easily. Thus $\lambda_{n}\rightarrow0$ as
$n\rightarrow+\infty$. For $\forall\varphi\in D_{0}^{1,2}(\Omega)$,
$$
\langle G'(v_{n}),\varphi\rangle=2a(v_{n},\varphi)+2\int_{\Omega}V(x) v_{n}\varphi dx+4b\int_{\Omega}|\nabla v_{n}|^{2}dx\int_{\Omega}\nabla v_{n}\nabla\varphi dx-6\int_{\Omega}v_{n}^{5}\varphi dx.
$$
Thus, by the boundedness of $|V|_{\frac{3}{2}}$ and the H\"{o}lder inequality, we have
$$\begin{array}{rcl}
|\langle G'(v_{n}),\varphi\rangle|&\leq&C_{1}\|v_{n}\|\|\varphi\|+4b\|v_{n}\|^{3}\|\varphi\|+6|v_{n}|_{6}^{5}|\varphi|_{6}\\
\\
&\leq&C_{1}\|v_{n}\|\|\varphi\|+4b\|v_{n}\|^{3}\|\varphi\|+C_{2}\|v_{n}\|^{5}\|\varphi\|.\\
\end{array}
$$
Then it follows from the boundedness of $\{v_{n}\}$ that $G'(v_{n})$
is bounded. Hence, $\lambda_{n}G'(v_{n})\rightarrow0$ and then
$I'(v_{n})\rightarrow0$ by \eqref{eq:vn-minimizing sequence of N}.
For $\forall\varphi\in D_{0}^{1,2}(\Omega)$, it is easy to see that
$$
\langle I'(u_{n})-I'(v_{n}),\varphi\rangle\rightarrow0\ \text{as}\  n\rightarrow+\infty.
$$
Thus, $I'(u_{n})\rightarrow0$.

By Lemma 3.1, there exists a number $k\in\mathbb{N}$ and a
subsequence of $\{u_n\}$, still denoted by $\{u_{n}\}$, such that
(\ref{eq:u(n)-decomposition})-(\ref{eq:I(u)-decomposition}) hold.
Then we have the following possibilities.\\
(1) If $u^{(0)}\not \equiv 0$ and  $k\geq 1$, by
\eqref{eq:u(0)-function}-(\ref{e2.2}), we have
\begin{equation}\label{eq:Psi'(u^(j))}
\langle I'(u^{(0)}),u^{(0)}\rangle\leq0,\qquad
\langle\Psi'(u^{(j)}),u^{(j)}\rangle\leq0,\quad 1\leq j\leq k.
\end{equation}
Now we show that there exists $t_{j}\in(0,1]$ ($0\leq j\leq k$) such that
\begin{equation}\label{eq:Psi'(u^(j)(x/t))}
\Big\langle
I'(u^{(0)}(\frac{\cdot}{t_{0}})),u^{(0)}(\frac{\cdot}{t_{0}})\Big\rangle=0,\quad\Big\langle\Psi'(u^{(j)}(\frac{\cdot}{t_{j}})),u^{(j)}(\frac{\cdot}{t_{j}})\Big\rangle=0,\quad
1\leq j\leq k,
\end{equation}
which implies that $I(u^{(0)}(\frac{\cdot}{t_{0}}))\geq M_\Omega=M_{\infty}$ and $\Psi(u^{(j)}(\cdot/t_{j}))\geq M_{\infty}$ for any $1\leq j\leq k$.
Now we prove (\ref{eq:Psi'(u^(j)(x/t))}). Let
$$
M(t):=\Big\langle
I'(u^{(0)}(\frac{\cdot}{t})),u^{(0)}(\frac{\cdot}{t})\Big\rangle,\quad
t\in(0,1].
$$
Then by direct calculation, we have
\begin{equation*}
M(t)=at\int_{\Omega}|\nabla
u^{(0)}|^2dx+t^3\int_{\Omega}V(tx)(u^{(0)})^2dx+bt^2\left(\int_{\Omega}|\nabla
u^{(0)}|^2dx\right)^2-t^3\int_{\Omega}(u^{(0)})^{6}dx.
\end{equation*}
Clearly, $M(t)$ is continuous in $(0,1]$, $M(t)>0$ for $t>0$ small enough, and $M(1)\leq 0$ by \eqref{eq:Psi'(u^(j))}, thus there exists
$t_{0}\in(0,1]$ such that $M(t_{0})=0$; that is,
$$
\Big\langle I'(u^{(0)}(\frac{\cdot}{t_{0}})),u^{(0)}(\frac{\cdot}{t_{0}})\Big\rangle=0.
$$
Similarly, there exists $t_{j}\in(0,1]$ ($1\leq j\leq k$) such that
$$
\Big\langle\Psi'(u^{(j)}(\frac{\cdot}{t_{j}})),u^{(j)}(\frac{\cdot}{t_{j}})\Big\rangle=0,\quad 1\leq j\leq k,
$$
thus \eqref{eq:Psi'(u^(j)(x/t))} is proved.
Moreover, by
\eqref{eq:u(0)-function} and \eqref{eq:u(i)-function}, we have
\begin{equation}\label{eq:u(0)-function-imply}
D_{1}:=\int_{\Omega}a|\nabla
u^{(0)}|^2+V(x)(u^{(0)})^2dx+bA^{2}\int_{\Omega}|\nabla
u^{(0)}|^2dx-\int_{\Omega}(u^{(0)})^{6}dx=0
\end{equation}
and
\begin{equation}\label{eq:u(i)-function-imply}
D_{2}:=\int_{\mathbb{R}^{3}}a|\nabla
u^{(i)}|^2dx+bA^{2}\int_{\mathbb{R}^{3}}|\nabla
u^{(i)}|^2dx-\int_{\mathbb{R}^{3}}(u^{(i)})^{6}dx=0
\end{equation}
with $1\leq i\leq k$. Then by \eqref{eq:Psi'(u^(j)(x/t))} and \eqref{eq:u(0)-function-imply}, we
have
\begin{equation}\begin{array}{rcl}\label{eq:widetilde I(u(0))-imply}
\widetilde{I}(u^{(0)})&=&\widetilde{I}(u^{(0)})-\mathlarger{\frac{1}{6}}D_{1}\\
\\
&=&\mathlarger{\frac{1}{3}\int_{\Omega}\left(a|\nabla u^{(0)}|^{2}+V(x)(u^{(0)})^{2}\right)dx+\frac{b}{12}A^2\int_{\Omega}|\nabla u^{(0)}|^{2}dx}\\
\\
&\geq&\mathlarger{\frac{t_{0}}{3}\int_{\Omega}a|\nabla u^{(0)}|^{2}dx+\frac{t_{0}^{3}}{3}\int_{\Omega}V(x) (u^{(0)})^{2}dx+\frac{bt_{0}^{2}}{12}\left(\int_{\Omega}|\nabla u^{(0)}|^{2}dx\right)^{2}}\\
\\
&=&I(u^{(0)}(\cdot/t_{0}))-\mathlarger{\frac{1}{6}}\big\langle I'(u^{(0)}(\cdot/t_{0})),u^{(0)}(\cdot/t_{0})\big\rangle\\
\\
&=&I(u^{(0)}(\cdot/t_{0}))\geq M_{\infty},
\end{array}
\end{equation}
and similarly, by \eqref{eq:Psi'(u^(j)(x/t))} and \eqref{eq:u(i)-function-imply}, we
have
\begin{equation}\label{eq:widetilde Psi(u(i))-imply}
\widetilde{\Psi}(u^{(i)})\geq \Psi(u^{(i)}(\cdot/t_{i}))\geq
M_{\infty}, \quad 1\leq i\leq k.
\end{equation}
Thus it follows from \eqref{eq:I(u)-decomposition},
\eqref{eq:widetilde I(u(0))-imply} and \eqref{eq:widetilde
Psi(u(i))-imply} that
$$
M_{\Omega}=\widetilde{I}(u^{(0)})+\sum_{i=1}^{k}\widetilde{\Psi}(u^{(i)})\geq(k+1)M_{\infty}\geq 2M_{\infty},
$$
which contradicts $M_{\Omega}=M_{\infty}$.\\
(2) If $u^{(0)}\equiv 0$ and  $k\geq 2$, similarly to the above
step,
$$
M_{\Omega}=\sum_{i=1}^{k}\widetilde{\Psi}(u^{(i)})\geq
kM_{\infty}\geq 2M_{\infty},
$$
which contradicts $M_{\Omega}=M_{\infty}$.\\
(3) If $u^{(0)}\not \equiv 0$ and  $k=0$, then $u^{(0)}$ is a ground
state solution of \eqref{eq:Omega}, which contradicts Lemma 2.4.
Thus, we must have
$$
u^{(0)}=0, \quad k=1, \quad A^{2}=|\nabla
u^{(1)}|^{2}_{L^{2}(\mathbb{R}^{3})},\quad
u_{n}(x)=u_{n}^{(0)}(x)+(\sigma_{n}^{1})^{-\frac{1}{2}}u_n^{(1)}\left(\frac{\cdot-y_{n}^{1}}{\sigma_{n}^{1}}\right),
$$
and $u^{(1)}$ satisfies
$$
-\left(a+b\int_{\mathbb{R}^{3}}|\nabla u^{(1)}|^{2}dx\right)\Delta
u^{(1)}=(u^{(1)})^5, \quad \text{in} \ \mathbb{R}^{3}
$$
with $\Psi(u^{(1)})=M_{\Omega}=M_{\infty}$; that is,
$u^{(1)}$ is a ground state solution of \eqref{eq:limit}. By Lemma 2.2,
$u^{(1)}$ does not change sign. By the positivity of $\{u_n\}$ and
$u_{n}^{(0)}\to u^{(0)}=0$  in $D_{0}^{1,2}(\Omega)$, we obtain
that $u^{(1)}\geq 0$. The strong maximum principle implies
$u^{(1)}>0$. By the above arguments and  Lemma 2.2, we obtain
(\ref{e3.1}). The proof is
complete. $\quad\Box$
\\
\textbf{Lemma 3.4.} {\it Let $a> 0$, $b\geq 0$,  $\Omega$ be as in (\ref{eq:Omega}), $V\in L^{3/2}(\Omega)$ be non-negative and continuous, and  $\{u_{n}\}\subset D_{0}^{1,2}(\Omega)$
be a positive  Palais-Smale sequence for $I$, i.e.,
$$
I(u_{n})\rightarrow c\  \ \text{and}\ \  I'(u_{n})\rightarrow 0,\  \ \text{as}\ \  n\rightarrow+\infty.
$$
Then if $c\in\big(M_{\infty},2M_{\infty}\big)$, $\{u_{n}\}$ is relatively compact.}\\
\textbf{Proof.} By Lemma 3.1, there exists a number $k\in\mathbb{N}$
and a  subsequence of $\{u_n\}$, still denoted by $\{u_{n}\}$, such
that \eqref{eq:u(n)-decomposition}-(\ref{eq:I(u)-decomposition})
hold. Similarly to the proof of Lemma 3.3, the cases $u^{(0)}\neq0$
and $k\geq1$ or $u^{(0)}=0$ and $k\geq2$ do not occur for
$c\in\big(M_{\infty},2M_{\infty}\big)$. If
$u^{(0)}=0$ and  $k=1$, by Lemma 2.2 and the discussion in the proof of Lemma 3.3, we have $c=M_{\infty}$, which contradicts
$c>M_{\infty}$.
In conclusion, we must have $u^{(0)}\neq0$ and $k=0$, which implies that $\{u_{n}\}$ is relatively compact. The lemma is proved. $\quad\Box$

\section{Proof of Theorem 1.1}
Let $\overline{u}_{\delta,y}$ given by \eqref{overline{u}_delta,y}  be the ground state solution of the equation \eqref{eq:limit}.
Without any loss of generality, we assume that $0\in\mathbb{R}^{3}\setminus\Omega$, then  $\overline{\rho}:=\inf\left\{\rho:\mathbb{R}^{3}\setminus\Omega\subset\overline{B_{\rho}(0)}\right\}>0$. We define
$$
v_{\bar{\rho}}(x):=\zeta(x)\overline{u}_{\delta,y}=\widetilde{\zeta}\left(\frac{|x|}{\bar{\rho}}\right)\overline{u}_{\delta,y},
$$
where  $\zeta$, $\widetilde{\zeta}$ are 
chosen as in the proof of Lemma 2.4.  
Define the operator $\Phi_{\bar{\rho}}:\mathbb{R}^{3}\times\mathbb{R}^{+}\rightarrow D^{1,2}(\mathbb{R}^{3})$ as follows
$$
\Phi_{\bar{\rho}}(y,\delta)=t_{v}v_{\bar{\rho}}(x),
$$
where $t_{v}$ is chosen such that $\big\langle
I'(\Phi_{\bar{\rho}}(y,\delta)),\Phi_{\bar{\rho}}(y,\delta)\big\rangle=0$.
 It is easy to see  that $\Phi_{\bar{\rho}}(y,\delta)$ and
$v_{\bar{\rho}}(x)$ vanish outside of $\Omega$, and thus they can be seen as elements in $D_{0}^{1,2}(\Omega)$ ($L^{6}(\Omega)$), and there hold
$$
\|\Phi_{\bar{\rho}}(y,\delta)\|=\|\Phi_{\bar{\rho}}(y,\delta)\|_{\mathbb{R}^{3}},\quad
\|v_{\bar{\rho}}\|=\|v_{\bar{\rho}}\|_{\mathbb{R}^{3}},
$$
$$
|\Psi_{\bar{\rho}}(y,\delta)|_{6}=|\Psi_{\bar{\rho}}(y,\delta)|_{L^{6}(\mathbb{R}^{3})},\quad |v_{\bar{\rho}}|_{6}=|v_{\bar{\rho}}|_{L^{6}(\mathbb{R}^{3})}.
$$
\textbf{Lemma 4.1.} {\it Let $a>0$, $b\geq0$,  $V\in L^{3/2}(\Omega)$ be non-negative and continuous. If $|V|_{\frac{3}{2}}\neq0$, then the functions $v_{\bar{\rho}}$ and $\Phi_{\bar{\rho}}(y,\delta)$ satisfy:
\begin{itemize}
  \item[$(1)$] $\Phi_{\bar{\rho}}(y,\delta)$ is continuous in $(y,\delta)$ for every $\bar{\rho}$;
  \item[$(2)$]  $I(v_{\bar{\rho}})\rightarrow M_{\infty}$, $\left\langle I'(v_{\bar{\rho}}),v_{\bar{\rho}}\right\rangle\to 0$, and $t_{v}\to 1$ as $|y|\rightarrow+\infty$, uniformly for every bounded $\bar{\rho}$, and bounded $\delta$ away from 0;
  \item[$(3)$] as $\bar{\rho}\rightarrow0$,  $I(v_{\bar{\rho}})\rightarrow M_{\infty}$, $\left\langle I'(v_{\bar{\rho}}),v_{\bar{\rho}}\right\rangle\to 0$, and $t_{v}\to 1$ as $\delta\to 0$ or $\delta\to+\infty$, uniformly in $y\in\mathbb{R}^{3}$.
\end{itemize}}
\noindent\textbf{Proof.} By the definition of $v_{\bar{\rho}}$ and (1) of Lemma 2.1, we know that (1) is right. $(2)$ follows from the same
arguments as in the proof of Lemma 2.4. Now  we  prove $(3)$. By
direct calculation, we have
$$\begin{array}{rcl}
\big\|v_{\bar{\rho}}-\overline{u}_{\delta,y}\big\|^{2}_{\mathbb{R}^{3}}&=&\mathlarger{\int}_{\mathbb{R}^{3}}\left|\nabla(v_{\bar{\rho}}-\overline{u}_{\delta,y})\right|^{2}dx\\
\\
&\leq&c_{1}\mathlarger{\int}_{B_{2\bar{\rho}}}\big|\nabla\overline{u}_{\delta,y}\big|^{2}dx     +c_{2}\mathlarger{\int}_{\bar{\rho}\leq|x|\leq2\bar{\rho}}\big|\overline{u}_{\delta,y}\nabla\zeta(x)\big|^{2}dx\\
\\
&\leq&c_{1}\mathlarger{\int}_{B_{2\bar{\rho}}}\big|\nabla\overline{u}_{\delta,y}\big|^{2}dx +\mathlarger{\frac{c_{3}}{\bar{\rho}^{2}}}\mathlarger{\int}_{\bar{\rho}\leq|x|\leq2\bar{\rho}}\big|\overline{u}_{\delta,y}\big|^{2}dx \\
\\
&\leq&c_{1}\mathlarger{\int}_{B_{2\bar{\rho}}}\big|\nabla\overline{u}_{\delta,y}\big|^{2}dx
+\mathlarger{\frac{c_{4}}{\bar{\rho}^{2}}}\left(\mathlarger{\int}_{\bar{\rho}\leq|x|\leq2\bar{\rho}}dx\right)^{2/3}
\left(\mathlarger{\int}_{\bar{\rho}\leq|x|\leq2\bar{\rho}}\big|\overline{u}_{\delta,y}\big|^{6}dx\right)^{1/3}\\
\\
&\rightarrow&0\quad \text{as}\ \ \bar{\rho}\rightarrow0,\ \ \text{for} \ \ \forall(y,\delta)\in\mathbb{R}^{3}\times\mathbb{R}^{+},
\end{array}
$$
where $c_{i}>0$ ($i=1,2,3,4$) are all constants; that is,
\begin{equation}\label{eq:|v_bar rho|-norm}
\big\|v_{\bar{\rho}}\big\|^{2}=\big\|v_{\bar{\rho}}\big\|^{2}_{\mathbb{R}^{3}}\rightarrow\big\|\overline{u}_{\delta,y}\big\|^{2}_{\mathbb{R}^{3}}\ \text{as}\ \bar{\rho}\rightarrow0, \quad \text{for}\ \forall(y,\delta)\in\mathbb{R}^{3}\times\mathbb{R}^{+}.
\end{equation}
As $\bar{\rho}\rightarrow0$, we have
\begin{equation}\label{eq:|v_bar rho| 6}
\begin{array}{rcl}
\left|\mathlarger{\int}_{\mathbb{R}^{3}}|v_{\bar{\rho}}|^{6}dx-\mathlarger{\int}_{\mathbb{R}^{3}}|\overline{u}_{\delta,y}|^{6}dx\right|&\le 
&\mathlarger{\int}_{B_{2\bar{\rho}}}\left| |v_{\bar{\rho}}|^{6}-|\overline{u}_{\delta,y}|^{6}\right|dx\\
\\
&\leq&c_{1}\mathlarger{\int_{B_{2\bar{\rho}}}\big|\overline{u}_{\delta,y}\big|^{6}dx}
\rightarrow 0,\quad \text{for}\  \forall(y,\delta)\in\mathbb{R}^{3}\times\mathbb{R}^{+};
\end{array}
\end{equation}
that is,
\begin{equation}\label{eq:|v_bar rho| 6}
\big|v_{\bar{\rho}}\big|_{6}^{6}=\big|v_{\bar{\rho}}\big|_{L^{6}(\mathbb{R}^{3})}^{6}\rightarrow\big|\overline{u}_{\delta,y}\big|_{L^{6}(\mathbb{R}^{3})}^{6}\ \ \text{as}\ \  \bar{\rho}\rightarrow0, \quad \text{for}\ \forall(y,\delta)\in\mathbb{R}^{3}\times\mathbb{R}^{+}.
\end{equation}
Similarly, by $V\in L^{3/2}(\Omega)$, (\ref{eq:|v_bar rho| 6}) and the H\"{o}lder inequality, we obtain that
\begin{equation}\label{eq:|V(x)v_bar rho| 2}
\int_{\Omega}V(x)|v_{\bar{\rho}}|^{2}dx\rightarrow\int_{\mathbb{R}^{3}}V(x)|\overline{u}_{\delta,y}|^{2}dx \ \ \text{as} \ \ \bar{\rho}\rightarrow0, \quad \text{for}\ \forall(y,\delta)\in\mathbb{R}^{3}\times\mathbb{R}^{+}.
\end{equation}
By $V\in L^{3/2}(\Omega)$, the definition of $\overline u_{\delta,y}$, and  arguing as in the proof of  \cite[Lemma 4.4]{Xie-1}, we know that for any $\varepsilon>0$, there exist
$\bar{\delta}_{1}=\bar{\delta}_{1}(\varepsilon)$ and $\bar{\delta}_{2}=\bar{\delta}_{2}(\varepsilon)$ such that
\begin{equation}\label{eq:|V(x)varphi|_2<varepsilon}
\int_{\mathbb{R}^{3}}V(x)|\overline{u}_{\delta,y}|^{2}dx=K\int_{\mathbb{R}^{3}}V(x)\varphi^{2}_{\delta,y}dx<\varepsilon,
\end{equation}
for $\forall y\in\mathbb{R}^{3}$ and  $\delta\in(0,\bar{\delta}_{1}]\cup[\bar{\delta}_{2},+\infty)$,
which combined with \eqref{eq:|V(x)v_bar rho| 2} gives that as $\bar{\rho}\rightarrow0$,
\begin{equation}\label{eq:|V(x)v_bar rho|<rightarrow0}
\int_{\Omega}V(x)|v_{\bar{\rho}}|^{2}dx\rightarrow0\ \ \text{as}\  \ \delta\rightarrow0 \  \ \text{or} \ \ \delta\rightarrow\infty,\quad \text{for} \ \forall y\in\mathbb{R}^{3}.
\end{equation}
By \eqref{eq:|v_bar rho|-norm}, \eqref{eq:|v_bar rho| 6}, \eqref{eq:|V(x)v_bar rho|<rightarrow0}, $\Psi(\overline{u}_{\delta,y})=M_{\infty}$, and $\langle\Psi'(\overline{u}_{\delta,y}),\overline{u}_{\delta,y}\rangle=0$, we obtain that as $\bar{\rho}\rightarrow0$, $I(v_{\bar{\rho}})\to M_{\infty}$ and
$\langle I'(v_{\bar{\rho}}),v_{\bar{\rho}}\rangle\to 0$ as $\delta\rightarrow0$ or $\delta\rightarrow\infty$ for $\forall y\in\mathbb{R}^{3}$.
Then similarly to the proof of \eqref{tn rightarrow1} in Lemma 2.4, we can show that as $\bar{\rho}\rightarrow0$, $t_{v}\rightarrow 1$ as $\delta\rightarrow0$ or $\delta\rightarrow\infty$ for $\forall y\in\mathbb{R}^{3}$. Thus, (3) is proved. $\quad\Box$\\
\textbf{Lemma 4.2.} {\it Let $a>0$, $b\geq0$, and $(V_{1})$ hold, then  there exists a small constant
$\widetilde{\rho}\in (0,1/8)$ such that for $\forall\bar{\rho}<\widetilde{\rho}$,}
$$
\sup_{(y,\delta)\in\mathbb{R}^{3}\times\mathbb{R}^{+}}I(\Phi_{\bar{\rho}}(y,\delta))<2M_{\infty}.
$$
\textbf{Proof.} By \eqref{eq:|v_bar rho|-norm}, \eqref{eq:|v_bar rho| 6}, \eqref{eq:|V(x)v_bar rho| 2}, the H\"{o}lder inequality, and
$$
\big\langle I'(\Phi_{\bar{\rho}}(y,\delta)),\Phi_{\bar{\rho}}(y,\delta)\big\rangle=\big\langle I'(t_{v}v_{\bar{\rho}}),t_{v}v_{\bar{\rho}}\big\rangle=0,
$$
we have
$$
\begin{array}{rcl}
t^{2}_{v}&=&\mathlarger{\frac{b\|v_{\bar{\rho}}\|^{4}+\sqrt{b^{2}\|v_{\bar{\rho}}\|^{8}+4|v_{\bar{\rho}}|_{6}^{6}\left(a\|v_{\bar{\rho}}\|^{2}+\mathlarger{\int}_{\Omega}V(x)v_{\bar{\rho}}^{2}dx\right)}}{2|v_{\bar{\rho}}|_{6}^{6}}}\\
\\
&\leq&\mathlarger{\frac{b\|v_{\bar{\rho}}\|^{4}+\sqrt{b^{2}\|v_{\bar{\rho}}\|^{8}+4|v_{\bar{\rho}}|_{6}^{6}\left(a\|v_{\bar{\rho}}\|^{2}+|V|_{\frac{3}{2}}|v_{\bar{\rho}}|_{6}^{2}\right)}}{2|v_{\bar{\rho}}|_{6}^{6}}}\\
\\
&\rightarrow&\mathlarger{\frac{b\|\overline{u}_{\delta,y}\|_{\mathbb{R}^{3}}^{4}+\sqrt{b^{2}\|\overline{u}_{\delta,y}\|_{\mathbb{R}^{3}}^{8}+4|\overline{u}_{\delta,y}|_{L^{6}(\mathbb{R}^{3})}^{6}\left(a\|\overline{u}_{\delta,y}\|_{\mathbb{R}^{3}}^{2}+|V|_{\frac{3}{2}}|\overline{u}_{\delta,y}|_{L^{6}(\mathbb{R}^{3})}^{2}\right)}}{2|\overline{u}_{\delta,y}|_{L^{6}(\mathbb{R}^{3})}^{6}}}
\end{array}
$$
as $\bar{\rho}\to 0$ for $\forall(y,\delta)\in\mathbb{R}^{3}\times\mathbb{R}^{+}$.
By (\ref{e3.2}), (\ref{overline{u}_delta,y}), (\ref{eq:|v_bar rho| 6}) and $\big\langle I'(\Phi_{\bar{\rho}}(y,\delta)),\Phi_{\bar{\rho}}(y,\delta)\big\rangle=0$, we have
$$
\begin{array}{rcl}
I(\Phi_{\bar{\rho}}(y,\delta))&=&\mathlarger{\frac{at^{2}_{v}}{2}\|v_{\bar{\rho}}\|^{2}+\frac{t_{v}^{2}}{2}\int_{\Omega}V(x)v_{\bar{\rho}}^{2}dx+\frac{bt^{4}_{v}}{4}\|v_{\bar{\rho}}\|^{4}-\frac{t^{6}_{v}}{6}|v_{\bar{\rho}}|_{6}^{6}}\\
\\
&=&\mathlarger{\frac{at^{2}_{v}}{3}\|v_{\bar{\rho}}\|^{2}+\frac{t_{v}^{2}}{3}\int_{\Omega}V(x)v_{\bar{\rho}}^{2}dx+\frac{bt^{4}_{v}}{12}\|v_{\bar{\rho}}\|^{4}}\\
\\
&<&\mathlarger{\frac{1}{3}\left(at^{2}_{v}\|v_{\bar{\rho}}\|^{2}+t_{v}^{2}\int_{\Omega}V(x)v_{\bar{\rho}}^{2}dx+bt^{4}_{v}\|v_{\bar{\rho}}\|^{4}\right)}=\mathlarger{\frac{1}{3}t^{6}_{v}|v_{\bar{\rho}}|_{6}^{6}},
\end{array}
$$
as $\bar{\rho}\to 0$ for any $(y,\delta)\in \mathbb{R}^3\times \mathbb{R}^+$. 
Let
$$
B=b\|\overline{u}_{\delta,y}\|_{\mathbb{R}^{3}}^{4}+\sqrt{b^{2}\|\overline{u}_{\delta,y}\|_{\mathbb{R}^{3}}^{8}
+4|\overline{u}_{\delta,y}|_{L^{6}(\mathbb{R}^{3})}^{6}\left(a\|\overline{u}_{\delta,y}\|_{\mathbb{R}^{3}}^{2}
+|V|_{\frac{3}{2}}|\overline{u}_{\delta,y}|_{L^{6}(\mathbb{R}^{3})}^{2}\right)}.
$$
Since $\mathlarger{|v_{\bar{\rho}}|_{6}^{6}}\rightarrow\mathlarger{|\overline{u}_{\delta,y}|_{L^{6}(\mathbb{R}^{3})}^{6}}$ as $\bar{\rho}\to 0$ for any $(y,\delta)\in \mathbb{R}^3\times \mathbb{R}^+$, and 
$$
\begin{array}{rcl}
\mathlarger{\frac{1}{3}t^{6}_{v}|v_{\bar{\rho}}|_{6}^{6}}&=&\mathlarger{\frac{1}{3}|\overline{u}_{\delta,y}|_{L^{6}(\mathbb{R}^{3})}^{6}}\cdot
\mathlarger{\left(\frac{B}{2|\overline{u}_{\delta,y}|_{L^{6}(\mathbb{R}^{3})}^{6}} \right)^3}
\\ \\
&=&\mathlarger{\frac{1}{3}\left(\frac{B}{2|\overline{u}_{\delta,y}|_{L^{6}(\mathbb{R}^{3})}^{4}} \right)^3}
=\mathlarger{\frac{1}{3}\left(\frac{b}{2}S^{2}+\sqrt{\frac{b}{4}S^{4}+aS+|V|_{\frac{3}{2}}}\right)^{3}},
\end{array}
$$
Under the condition $(V_{1})$, we then obtain  that
$$
\sup_{(y,\delta)\in\mathbb{R}^{3}\times\mathbb{R}^{+}}I(\Phi_{\bar{\rho}}(y,\delta))<2M_{\infty}\ \text{as} \ \ \bar{\rho}\to 0 \ \ \text{for} \ \ \forall(y,\delta)\in\mathbb{R}^{3}\times\mathbb{R}^{+}.
$$
The lemma is proved. $\quad\Box$

From now on, we will assume $\Omega$ fixed in  such a way that
$$
{\rm diam}(\mathbb{R}^{3}\setminus\Omega)<\widetilde{\rho},
$$
where ${\rm diam}(D)$ is defined by ${\rm diam}(D):=\sup\{|x-y|: \ x,y\in D\}$, $\widetilde{\rho}\in (0,1/8)$ is the constant obtained in  Lemma 4.2. Then for $\forall x_{0}\in\mathbb{R}^{3}\setminus\Omega$,
$$
\mathbb{R}^{3}\setminus\Omega\subset B_{\widetilde{\rho}}(x_{0}):=\{x\in\mathbb{R}^{3}:|x-x_{0}|<\widetilde{\rho}\}.
$$
Thus, we have 
$$
\overline{\rho}:=\inf\{\rho:\mathbb{R}^{3}\setminus\Omega\subset\overline{B_{\rho}(0)}\}<\widetilde{\rho}<1/8,
$$
and  $\mathbb{R}^{3}\setminus\Omega\subset B_{1/8}:=\{x\in\mathbb{R}^{3}:|x|<1/8\}$.
Define $\chi_{1}, \chi_{2}:\mathbb{R}^{+}\to \mathbb{R}$ as follows
$$
\chi_{1}(t):=\left\{
\begin{array}{ll}
4, \quad t\leq\mathlarger{\frac{1}{4}},\\
\\
1/t,\quad t>\mathlarger{\frac{1}{4}},\\
\end{array}\right.\ \text{and}\ \ 
\chi_{2}(t):=\left\{
\begin{array}{ll}
0, \quad t<1,\\
\\
1, \quad t\geq1.\\
\end{array}\right.
$$
Now we give a barycenter function, which is similar to that of \cite{Xie-1}, but in order to adapt it to our problem in exterior domains, some appropriate changes are needed. Define $\beta:D^{1,2}(\mathbb{R}^{3})\rightarrow\mathbb{R}^{3}\times\mathbb{R}^{+}$ as
$$
\beta(u):=\frac{1}{S^{\frac{3}{2}}K}\int_{\mathbb{R}^{3}}\left(\chi_{1}(|x|)x,\chi_{2}(|x|)\right)|\nabla u|^{2}dx=(\beta_{1}(u),\beta_{2}(u)),
$$
where
$$
\beta_{1}(u):=\frac{1}{S^{\frac{3}{2}}K}\int_{\mathbb{R}^{3}}\chi_{1}(|x|)x|\nabla u|^{2}dx,\quad
\beta_{2}(u):=\frac{1}{S^{\frac{3}{2}}K}\int_{\mathbb{R}^{3}}\chi_{2}(|x|)|\nabla u|^{2}dx,
$$
and $K=\frac{1}{2}\left(bS^{\frac{3}{2}}+\sqrt{b^{2}S^{3}+4a}\right)$. Let $\mathfrak{B}_{0}$ be the
subset of $D_{0}^{1,2}(\Omega)$ defined by
$$
\mathfrak{B}_{0}:=\left\{u\in \mathcal{N}:\beta(u)=(0,0,0,{1}/{2})\right\},
$$
where $\mathcal{N}$ is defined in Section 2.\\
\textbf{Lemma 4.3.} {\it If $|y|\geq\frac{1}{2}$, then}
$$
\beta_{1}(\overline{u}_{\delta,y})=\frac{y}{|y|}+o(1) \ \ as \ \ \delta\rightarrow0.
$$
\textbf{Proof.} Fix $|y|\geq\frac{1}{2}$. For any  $\varepsilon>0$ small enough, we have $B_{1/4}(0)\cap B_{\varepsilon}(y)=\emptyset$. By  the definition of $\overline u_{\delta,y}$, we have
\begin{equation}\label{eq:B_1/2(0)-o(1)}
\frac{1}{S^{\frac{3}{2}}K}\int_{B_{1/4}(0)}4x|\nabla \overline{u}_{\delta,y}|^{2}dx=\frac{1}{S^{\frac{3}{2}}}\int_{B_{1/4}(0)}4x|\nabla \varphi_{\delta,y}|^{2}dx\rightarrow0 \ \text{as}\ \delta\rightarrow0
\end{equation}
and
\begin{equation}\label{eq:D-o(1)}
\frac{1}{S^{\frac{3}{2}}K}\int_{D}\frac{x}{|x|}|\nabla \overline{u}_{\delta,y}|^{2}dx=\frac{1}{S^{\frac{3}{2}}}\int_{D}\frac{x}{|x|}|\nabla \varphi_{\delta,y}|^{2}dx\rightarrow0\ \text{as} \ \delta\rightarrow0,
\end{equation}
where $D=\mathbb{R}^{3}\backslash(B_{1/4}(0)\cup B_{\varepsilon}(y))$. And for any $x\in B_{\varepsilon}(y)$, considering $|y|\geq\frac{1}{2}$ and $\varepsilon>0$ small enough, we have
$$
\left|\frac{x}{|x|}-\frac{y}{|y|}\right|<c_{1}\varepsilon,
$$
where $c_{1}>0$ is a constant.  This implies that
\begin{equation}\label{eq:B_varepsilon(y)-o(1)}
\begin{array}{cl}
&\quad\mathlarger{\left|\frac{y}{|y|}-\frac{1}{S^{\frac{3}{2}}K}\int_{B_{\varepsilon}(y)}\frac{x}{|x|}|\nabla \overline{u}_{\delta,y}|^{2}dx\right|}\\
\\
&=\mathlarger{\left|\frac{1}{S^{\frac{3}{2}}K}\int_{B_{\varepsilon}(y)}\left(\frac{y}{|y|}-\frac{x}{|x|}\right)|\nabla \overline{u}_{\delta,y}|^{2}dx+\frac{1}{S^{\frac{3}{2}}K}\int_{B_{\varepsilon}^c(y)}\frac{y}{|y|}|\nabla \overline{u}_{\delta,y}|^{2}dx\right|}\\
\\
&\leq\mathlarger{\frac{c_{1}\varepsilon}{S^{\frac{3}{2}}K}\int_{\mathbb{R}^{3}}|\nabla \overline{u}_{\delta,y}|^{2}dx+     \frac{1}{S^{\frac{3}{2}}K}\int_{B_{\varepsilon}^c(y)}|\nabla \overline{u}_{\delta,y}|^{2}dx}
\rightarrow0\quad \text{as}\ \ \delta\rightarrow0.
\end{array}
\end{equation}
Therefore, by \eqref{eq:B_1/2(0)-o(1)}, \eqref{eq:D-o(1)} and \eqref{eq:B_varepsilon(y)-o(1)},
$$
\begin{array}{cl}
\beta_{1}(\overline{u}_{\delta,y})&=\mathlarger{\frac{1}{S^{\frac{3}{2}}K}\left[\int_{B_{1/4}(0)}4x|\nabla \overline{u}_{\delta,y}|^{2}dx+\int_{D}\frac{x}{|x|}|\nabla \overline{u}_{\delta,y}|^{2}dx+\int_{B_{\varepsilon}(y)}\frac{x}{|x|}|\nabla \overline{u}_{\delta,y}|^{2}dx\right]}\\
\\
&\rightarrow \mathlarger{\frac{y}{|y|}}\quad as \ \delta\rightarrow0.
\end{array}
$$
The lemma is proved. $\quad\Box$\\
\textbf{Lemma 4.4.} {\it Assume that $a>0$, $b\geq0$, $V\in L^{3/2}(\Omega)$ is non-negative and continuous, $|V|_{\frac{3}{2}}\neq0$, and $c_{0}:=\inf_{u\in
\mathfrak{B}_{0}\cap P}I(u)$, where $P$ is the cone of non-negative
functions of $D_{0}^{1,2}(\Omega)$. Then
\begin{equation}\label{eq:c0>ground state}
c_{0}>M_{\infty}
\end{equation}
and as $\bar{\rho}\rightarrow0$, there exist $R_{0}>1/2$, $\delta_{1}\in(0,1/2),\delta_{2}>1/2$ such that}
\begin{equation}\label{eq:Lemma 4.3-(a)(b)}
\begin{array}{rcl}
&(a)&\quad \mathlarger{\beta_{2}\left(\Phi_{\bar{\rho}}(y,\delta)\right)<\frac{1}{2},\quad  \ \forall|y|<\frac{1}{2}\ \text{and}\ \delta\leq\delta_{1}};\\
\\
&(b)&\quad \mathlarger{\left|\beta_{1}\left(\Phi_{\bar{\rho}}(y,\delta)\right)-\frac{y}{|y|}\right|<\frac{1}{4},\quad \ \forall|y|\geq\frac{1}{2}\ \text{and}\  \delta\leq\delta_{1}};\\
\\
&(c)&\quad \mathlarger{\beta_{2}\left(\Phi_{\bar{\rho}}(y,\delta)\right)>\frac{1}{2},\quad \forall y\in\mathbb{R}^{3}\ \text{and}\ \delta\geq\delta_{2}};\\
\\
&(d)&\quad I\big(\Phi_{\bar{\rho}}(y,\delta)\big)<(c_{0}+M_{\infty})/2,\quad \forall |y|\in \mathbb{R}^3,\ \delta=\delta_1\ \text{or}\ \delta= \delta_2;\\
\\
&(e)&\quad I\big(\Phi_{\bar{\rho}}(y,\delta)\big)\in\left(M_{\infty},[c_{0}+M_{\infty}]\big/2\right),\quad \forall |y|\geq R_{0} \ \text{and}\ \delta_1\leq \delta\leq \delta_2;\\
\\
&(f)&\quad \Big(\beta_{1}(\Phi_{\bar{\rho}}(y,\delta)),y\Big)_{\mathbb{R}^{3}}>0,\quad  \forall |y|= R_{0}\ \text{and}\ \delta\in[\delta_{1},\delta_{2}].\\
\end{array}
\end{equation}
\textbf{Proof.} Obviously, $c_{0}\geq M_{\infty}$. To
prove \eqref{eq:c0>ground state}, we suppose
$c_{0}=M_{\infty}$ by contradiction. Then there exists
a sequence $\{u_{n}\}\subset D_{0}^{1,2}(\Omega)\cap P$ satisfying
\begin{equation}\label{I(u_{n})-ps}
\langle I'(u_{n}),u_{n}\rangle=0,\quad\beta(u_{n})=(0,0,0,1/2),\
\text{for\ all}\  n
\end{equation}
and $I(u_{n})\rightarrow M_{\infty}$ as $n\to +\infty$.
Moreover, $u_{n}$ is not relatively compact because
$M_{\infty}$ is not achieved in $D_{0}^{1,2}(\Omega)$ by
Lemma 2.4. Then by Lemma 3.3, there exists a positive ground state
solution $\overline{u}_{\delta_{n},y_{n}}$ of \eqref{eq:limit} such that
$$
u_{n}(x)=\overline{u}_{\delta_{n},y_{n}}(x)+w_{n}(x),\quad\forall x\in\mathbb{R}^{3},
$$
where $\{y_{n}\}\subset\mathbb{R}^{3}$, $\{\delta_{n}\}\subset\mathbb{R}^{+}$, and $\{w_{n}\}$ is a sequence converging strongly to $0$  in $D^{1,2}(\mathbb{R}^{3})$.
For $\overline{u}_{\delta_{n},y_{n}}$, there exists a subsequence $(\delta_{n},y_{n})$ such that one of these cases occurs:
\begin{equation}\label{eq:(delta_{n},y_{n})}
\begin{array}{rcl}
&(i)&\quad \delta_{n}\rightarrow+\infty \ \ \text{as} \ \ n\rightarrow+\infty;\\
\\
&(ii)&\quad \delta_{n}\rightarrow\delta\neq0 \ \  \text{as} \ \ n\rightarrow+\infty;\\
\\
&(iii)&\quad \delta_{n}\rightarrow0\ \ \text{and}\ \ y_{n}\rightarrow y\ \ \text{with}\ \ |y|<1/2, \ \  \text{as} \ \ n\rightarrow+\infty;\\
\\
&(iv)&\quad \delta_{n}\rightarrow0 \ \ \text{as}\ \ n\rightarrow+\infty,\quad \text{and} \ \ |y_{n}|\geq 1/2 \ \  \text{for} \ \ n \ \ \text{large}.\\
\end{array}
\end{equation}
In fact, none of the possibilities \eqref{eq:(delta_{n},y_{n})} can be true. Because as $n\to+\infty$, we have that
$$
\beta(\overline{u}_{\delta_{n},y_{n}}+w_{n})=\beta(\overline{u}_{\delta_{n},y_{n}})+o_{n}(1).
$$
Then it follows from \eqref{I(u_{n})-ps} that
\begin{equation}\label{eq:beta(u{n})=(0,1/2)}
\beta_{1}(\overline{u}_{\delta_{n},y_{n}})\rightarrow (0,0,0)\  \ \text{and}\  \  \beta_{2}(\overline{u}_{\delta_{n},y_{n}})\rightarrow\frac{1}{2},\quad \text{as} \ n\rightarrow+\infty.
\end{equation}
If \eqref{eq:(delta_{n},y_{n})} $(i)$ holds, then
$$\begin{array}{rcl}
\beta_{2}(\overline{u}_{\delta_{n},y_{n}})&=&\mathlarger{\frac{1}{S^{\frac{3}{2}}K}\int_{\mathbb{R}^{3}\setminus B_{1}(0)}|\nabla\overline{u}_{\delta_{n},y_{n}}|^{2}dx}\\
\\
&=&\mathlarger{1-\frac{1}{S^{\frac{3}{2}}K}\int_{B_{1}(0)}|\nabla\overline{u}_{\delta_{n},y_{n}}|^{2}dx}
=1+o(1),\quad n\rightarrow+\infty,
\end{array}
$$
clearly, it contradicts \eqref{eq:beta(u{n})=(0,1/2)}. If \eqref{eq:(delta_{n},y_{n})} $(ii)$ holds, then $|y_{n}|\rightarrow+\infty$, because if not, $\overline{u}_{\delta_{n},y_{n}}$ would converge strongly in $D^{1,2}(\mathbb{R}^{3})$ which contradicts to Lemma 2.4. Then we have
$$\begin{array}{rcl}
\beta_{2}(\overline{u}_{\delta_{n},y_{n}})&=&\beta_{2}(\overline{u}_{\delta,y_{n}})+o(1)\\
\\
&=&\mathlarger{\frac{1}{S^{\frac{3}{2}}K}\int_{\mathbb{R}^{3}\setminus B_{1}(y_{n})}|\nabla\overline{u}_{\delta,0}|^{2}dx+o(1)}
\\ \\
&=&\mathlarger{1-\frac{1}{S^{\frac{3}{2}}K}\int_{B_{1}(y_{n})}|\nabla\overline{u}_{\delta,0}|^{2}dx+o(1)}
=1+o(1),\quad n\rightarrow+\infty,
\end{array}
$$
which contradicts \eqref{eq:beta(u{n})=(0,1/2)}. If \eqref{eq:(delta_{n},y_{n})} $(iii)$ holds, then
$$\begin{array}{rcl}
\beta_{2}(\overline{u}_{\delta_{n},y_{n}})&=&\mathlarger{\frac{1}{S^{\frac{3}{2}}K}\int_{\mathbb{R}^{3}\setminus B_{1}(0)}|\nabla\overline{u}_{\delta_{n},y_{n}}|^{2}dx}\\
\\
&=&\mathlarger{\frac{1}{S^{\frac{3}{2}}K}\int_{\mathbb{R}^{3}\setminus B_{1}(y_{n})}|\nabla\overline{u}_{\delta_{n},0}|^{2}dx}
=o(1),\quad n\rightarrow+\infty,
\end{array}
$$
which contradicts \eqref{eq:beta(u{n})=(0,1/2)}. If \eqref{eq:(delta_{n},y_{n})} $(iv)$ holds, then by Lemmas 4.3, we have
$$
\beta_{1}(\overline{u}_{\delta_{n},y_{n}})=\frac{y_{n}}{|y_{n}|}+o(1),\quad n\rightarrow+\infty,
$$
which contradicts \eqref{eq:beta(u{n})=(0,1/2)}. Thus \eqref{eq:c0>ground state} is proved.

Nextly, we divide the proof of \eqref{eq:Lemma 4.3-(a)(b)} into several steps.

(1) By Lemma 4.1 (3), as $\bar{\rho}\rightarrow0$, we have
$$
I(\Phi_{\bar{\rho}}(y,\delta))=I(t_{v}v_{\bar{\rho}})\rightarrow M_{\infty} \ \text{as}\ \delta \to 0,\  \text{uniformly\ in\ } y\in \mathbb{R}^3.
$$
Then by $\Phi_{\bar{\rho}}(y,\delta)\in \mathcal{N}$ and Lemma 3.3, as $\bar{\rho}\rightarrow0$, we have
\begin{equation}\label{eq:beta_2(Phi(y,delta))}
\beta_{2}(\Phi_{\bar{\rho}}(y,\delta))\rightarrow\beta_{2}(\overline{u}_{\delta,y}(x))\ \text{as} \ \delta\rightarrow0, \  \text{uniformly\ in\ } y\in\mathbb{R}^3.
\end{equation}
By  the definition of $\overline u_{\delta,y}$,  we have
$$
\begin{array}{rcl}
\beta_{2}(\overline{u}_{\delta,y}(x))&=&\mathlarger{\frac{1}{S^{\frac{3}{2}}K}\int_{\mathbb{R}^{3}\backslash B_{1}(0)}|\nabla\overline{u}_{\delta,y}(x)|^{2}dx}\\
\\
&=&\mathlarger{\frac{1}{S^{\frac{3}{2}}}\int_{\mathbb{R}^{3}\backslash B_{1}(y)}|\nabla\varphi_{\delta,0}(x)|^{2}dx}
=o(1)\ \text{as}\ \delta\to 0, \ \text{for\ }\forall |y|<1/2.
\end{array}
$$
Thus, as $\bar{\rho}\to0$, there exists a small $\delta_{1}\in(0,1/2)$ such that \eqref{eq:Lemma 4.3-(a)(b)} (a) holds.

(2) Similarly to (1), by  Lemma 4.1 (3) and Lemma 3.3, as $\bar{\rho}\to0$, we have
$$
\beta_{1}(\Phi_{\bar{\rho}}(y,\delta))\rightarrow\beta_{1}(\overline{u}_{\delta,y}(x))\ \text{as} \ \delta\rightarrow0, \text{\ uniformly\ in\ }y\in\mathbb{R}^3.
$$
Then from Lemma 4.3, as $\bar{\rho}\to0$, there exists a small $\delta_{1}\in (0,1/2)$ such that \eqref{eq:Lemma 4.3-(a)(b)} (b) holds.

(3) Similarly to (1), by  Lemma 4.1 (3) and Lemma 3.3, as $\bar{\rho}\to0$, we have
\begin{equation}\label{e1.3}
\beta_{2}(\Phi_{\bar{\rho}}(y,\delta))\rightarrow\beta_{2}(\overline{u}_{\delta,y}(x))\  \ \text{as} \ \delta\rightarrow+\infty, \  \ \text{uniformly\ in\ } y\in\mathbb{R}^3.
\end{equation}
By \eqref{eq:varphi_delta,x0}, it is easy to see that
$$
\int_{B_{1}(0)}|\nabla\varphi_{\delta,y}|^{2}dx\to0\ \ \text{as} \ \ \delta\rightarrow+\infty,
$$
then
$$
\beta_{2}(\overline{u}_{\delta,y}(x))=1-\frac{1}{S^{\frac{3}{2}}}\int_{B_{1}(0)}|\nabla\varphi_{\delta,y}|^{2}dx=1+o(1)\  \text{as} \ \delta\rightarrow+\infty.
$$
Thus, as $\bar{\rho}\rightarrow0$, by (\ref{e1.3}), there exists a large $\delta_{2}>1/2$ such that
\eqref{eq:Lemma 4.3-(a)(b)} (c) holds.

(4)  By Lemma 4.1 (3), as $\bar{\rho}\to0$, we have
$$
I(\Phi_{\bar{\rho}}(y,\delta))=I(t_{v}v_{\bar{\rho}})\rightarrow M_{\infty}\  \ \text{as}\ \ \delta \to0\ \text{or}\  \ \delta\to+\infty,\ \text{uniformly\ in\ }y\in\mathbb{R}^3,
$$
which combined with \eqref{eq:c0>ground state} gives that as $\bar{\rho}\to0$, there exists $\delta_1\in (0,1/2)$ and $\delta_2>1/2$ such that
\eqref{eq:Lemma 4.3-(a)(b)} (d) holds.

(5) By $\Phi_{\bar{\rho}}(y,\delta)\in \mathcal{N}$, we have $I(\Phi_{\bar{\rho}}(y,\delta))\geq M_{\Omega}$, and then by Lemma 2.4, we have
$$
I(\Phi_{\bar{\rho}}(y,\delta))>M_{\Omega}=M_{\infty},\quad\forall y\in\mathbb{R}^{3}, \delta>0.
$$
By  Lemma 4.1 (2), we have
$$
I(\Phi_{\bar{\rho}}(y,\delta))=I(t_{v}v_{\bar{\rho}})\rightarrow M_{\infty}\ \text{as} \ |y|\rightarrow+\infty,
$$
uniformly for every bounded $\bar{\rho}$, and bounded $\delta$ away from 0. This and \eqref{eq:c0>ground state} imply that there exists  $R_0>1/2$ such that
$$
I(\Phi_{\bar{\rho}}(y,\delta))<\frac{c_{0}+M_{\infty}}{2},\quad \forall |y|\geq R_0, \ \ \delta\in[\delta_1,\delta_2], \ \ \bar{\rho}<1.
$$
Thus, \eqref{eq:Lemma 4.3-(a)(b)} (e) is satisfied.

(6) By $\Phi_{\bar{\rho}}(y,\delta)\in \mathcal{N}$, Lemma 4.1 (2), and Lemma 3.3, we have
\begin{equation}\label{e1.4}
\big(\beta_{1}(\Phi_{\bar{\rho}}(y,\delta)),y\big)_{\mathbb{R}^{3}}\rightarrow \big(\beta_{1}(\overline{u}_{\delta,y}),y\big)_{\mathbb{R}^{3}}\  \text{as} \ \ |y|\rightarrow+\infty,
\end{equation}
uniformly in $\delta\in[\delta_1,\delta_2]$ and $\bar{\rho}<1$. Set
$$
(\mathbb{R}^{3})_{y}^{+}:=\{x\in\mathbb{R}^{3}:(x,y)_{\mathbb{R}^{3}}>0\}\ \  \text{and}\  \ (\mathbb{R}^{3})_{y}^{-}:=\mathbb{R}^{3}\setminus\big(\mathbb{R}^{3}\big)_{y}^{+}.
$$
Since $\delta$ varies in the compact set $[\delta_{1},\delta_{2}]$,  there exist a large  $R_{0}\in\mathbb{R}^{+}$ and $\overline{r}\in(0,1/4)$ such that if $|y|\geq R_{0}$, the ball
$B_{\overline{r}}(\tilde{y})=\{x\in\mathbb{R}^{3}:|x-\tilde{y}|<\overline{r}\}\subset(\mathbb{R}^{3})_{y}^{+}$
with $\tilde{y}$ satisfying $|y-\tilde{y}|=\frac{1}{2}$ and by  the definition of $\overline u_{\delta,y}$,  we have
$$
|\nabla\varphi_{\delta,y}(x)|^{2}\geq C_{0}>0,\ \forall x\in B_{\overline{r}}(\tilde{y})
$$
and
$$
|\nabla\varphi_{\delta,y}(x)|^{2}\leq \frac{C_{1}}{|x-y|^{4}}, \ \forall x\in(\mathbb{R}^{3})_{y}^{-},
$$
where $C_{0},C_{1}$ are positive constants.
Hence, for any $|y|\geq R_{0}$, $\delta\in [\delta_1,\delta_2]$, we have
\begin{equation}\label{eq:(beta_{1}(u(x-y_{n})),y_{n})}
\begin{array}{cl}
\big(\beta_{1}(\overline{u}_{\delta,y}(x)),y\big)_{\mathbb{R}^{3}}
&=\mathlarger{\frac{1}{S^{\frac{3}{2}}}\int_{(\mathbb{R}^{3})_{y}^{+}}|\nabla\varphi_{\delta,y}(x)|^{2}\chi_{1}(|x|)(x,y)_{\mathbb{R}^{3}}dx}\\
\\
&\qquad\qquad+\mathlarger{\frac{1}{S^{\frac{3}{2}}}\int_{(\mathbb{R}^{3})_{y}^{-}}|\nabla\varphi_{\delta,y}(x)|^{2}\chi_{1}(|x|)(x,y)_{\mathbb{R}^{3}}dx}\\
\\
&\qquad\geq \mathlarger{C_2\int_{B_{\bar{r}}(\tilde{y})}\frac{1}{|x|}|x||y|dx-C_3\int_{(\mathbb{R}^{3})_{y}^{-}}\frac{1}{|x-y|^4}|y|dx}\\
\\
&\qquad\geq C_{4}|y|+C_5>0,
\end{array}
\end{equation}
where $C_{i}(i=2,3,4,5)$ are positive constants.
So we can choose $R_{0}>1/2$ such that  \eqref{eq:Lemma 4.3-(a)(b)} (f) holds.
The lemma is proved. $\quad\Box$

Now we define a bounded domain $D_{y,\delta}$ by
$$
D_{y,\delta}:=\Big\{(y,\delta)\in\mathbb{R}^{3}\times\mathbb{R}^{+}:|y|\leq R_{0},\ \delta\in[\delta_{1},\delta_{2}]\Big\},
$$
where $\delta_{1}$, $\delta_{2}$ and $R_{0}$ are constants found in Lemma 4.4. Note that $\partial D_{y,\delta}=D^{1}\cup D^{2}\cup D^{3}\cup D^{4}$, where
\begin{equation}\label{partial D(y,delta)}
\begin{array}{cl}
&D^{1}:=\left\{(y,\delta):|y|<\frac{1}{2},\ \delta=\delta_{1}\right\},\\
\\
&D^{2}:=\left\{(y,\delta):\frac{1}{2}\leq|y|\leq R_{0},\ \delta=\delta_{1}\right\},\\
\\
&D^{3}:=\left\{(y,\delta):|y|\leq R_{0},\ \delta=\delta_{2}\right\},\\
\\
&D^{4}:=\left\{(y,\delta):|y|= R_{0},\ \delta\in[\delta_{1},\delta_{2}]\right\}.\\
\end{array}
\end{equation}
We consider the subset $\Sigma$ of $D_{0}^{1,2}(\Omega)$ defined
by
$$
\Sigma:=\left\{\Phi_{\bar{\rho}}(y,\delta):\ (y,\delta)\in D_{y,\delta}\right\},
$$
then by the definition of $\Phi_{\bar{\rho}}(y,\delta)$, we have  $\Sigma\subset P\cap \mathcal{N}$, here $P$ is the cone of all non-negative
functions in  $D_{0}^{1,2}(\Omega)$. Define
$$
H:=\left\{h: h\in C(P\cap \mathcal{N}, P\cap \mathcal{N}), \
h(u)=u\ \text{for\ any\ } u \ \text{with\ }
I(u)<\frac{c_{0}+M_{\infty}}{2}\right\},
$$
and $\Gamma:=\{A\subset P\cap \mathcal{N}:A=h(\Sigma),\  h\in H\}$.\\
\textbf{Lemma 4.5.} {\it Assume that $a>0$, $b\geq0$, $V\in L^{3/2}(\Omega)$ is non-negative and continuous. Let $A\in\Gamma$, then as $\bar{\rho}\rightarrow0$,}
$$
A\cap \mathfrak{B}_{0}\neq\emptyset.
$$
\textbf{Proof.} It is equivalent to show that for any $ h\in H$, as $\bar{\rho}\rightarrow0$,
there exists $(\widetilde{y},\widetilde{\delta})\in D_{y,\delta}$ such that
\begin{equation}\label{eq:(beta h Phi)=(0,1/2)}
(\beta\circ h\circ\Phi_{\bar{\rho}})(\widetilde{y},\widetilde{\delta})=(0,0,0,{1}/{2}).
\end{equation}
With $h$ being chosen arbitrarily, we put $\varphi_{h}:\mathbb{R}^{3}\times\mathbb{R}^{+}\rightarrow\mathbb{R}^{3}\times\mathbb{R}^{+}$ with
$
\varphi_{h}:=\beta\circ h\circ\Phi_{\bar{\rho}},
$
and define $\varphi:D_{y,\delta}\rightarrow\mathbb{R}^{3}\times\mathbb{R}^{+}$ by
$
\varphi:=\beta\circ\Phi_{\bar{\rho}}.
$

Firstly, we claim that as $\bar{\rho}\rightarrow0$,
\begin{equation}\label{eq:d(varphi_h)=d(varphi)}
d\left(\varphi_{h},D_{y,\delta},(0,0,0,{1}/{2})\right)=d\left(\varphi,D_{y,\delta},(0,0,0,{1}/{2})\right).
\end{equation}
Indeed, by \eqref{eq:Lemma 4.3-(a)(b)}(d), \eqref{eq:Lemma 4.3-(a)(b)}(e) and \eqref{partial D(y,delta)}, we have as $\bar{\rho}\rightarrow0$,
$$
I(\Phi_{\bar{\rho}}(y,\delta))<\frac{c_{0}+M_{\infty}}{2},\quad\forall(y,\delta)\in \partial D_{y,\delta}.
$$
Thus, as $\bar{\rho}\rightarrow0$,
$$
h(\Phi_{\bar{\rho}}(y,\delta))=\Phi_{\bar{\rho}}(y,\delta), \quad\forall(y,\delta)\in \partial D_{y,\delta},
$$
and then
$$
\varphi_{h}(y,\delta)=(\beta\circ\Phi_{\bar{\rho}})(y,\delta)=\varphi(y,\delta),\quad \forall(y,\delta)\in \partial D_{y,\delta}.
$$
This proves \eqref{eq:d(varphi_h)=d(varphi)}.

Nextly, we prove that as $\bar{\rho}\rightarrow0$,
\begin{equation}\label{eq:d(varphi)=1}
d\left(\varphi,D_{y,\delta},\left(0,0,0,{1}/{2}\right)\right)=1.
\end{equation}
Consider the homotopy
$$
G(y,\delta,t)=t\varphi(y,\delta)+(1-t)(y,\delta),\quad 0\leq t\leq1,
$$
where $\varphi$ is homotopic to the identity $Id_{D_{y,\delta}}$.
By the homotopy invariance property of the topological degree and the fact that
$$
d\left(Id_{D_{y,\delta}},D_{y,\delta},\left(0,0,0,{1}/{2}\right)\right)=1,
$$
it suffices to prove that as $\bar{\rho}\rightarrow0$,
\begin{equation}\label{eq:varphi(y,delta)neq(0,1/2)}
G(y,\delta,t)\neq(0,0,0,{1}/{2}),\quad \text{for \ }\forall (y,\delta)\in \partial D_{y,\delta}\ \text{and \ } t\in[0,1].
\end{equation}
Now we prove (\ref{eq:varphi(y,delta)neq(0,1/2)}). By \eqref{eq:Lemma 4.3-(a)(b)} (a) and $\delta_1<1/2$, for  $\forall (y,\delta)\in D^{1}$, we have
\begin{equation}\label{eq:D^1-beta_2}
(1-t)\delta_{1}+t(\beta_{2}\circ\Phi_{\bar{\rho}})(y,\delta_{1})<\frac{1}{2},\quad \forall t\in[0,1].
\end{equation}
By \eqref{eq:Lemma 4.3-(a)(b)} (b) and $\delta_1<1/2$, for $\forall (y,\delta)\in D^{2}$, we have
$$
\left|\beta_{1}(\Phi_{\bar{\rho}}(y,\delta))-\frac{y}{|y|}\right|<\frac{1}{4}.
$$
Thus,
\begin{equation}\label{eq:D^2-beta_1}
\begin{array}{cl}
\left|(1-t)y+t(\beta_{1}\circ\Phi_{\bar{\rho}})(y,\delta_{1})\right|&\geq\mathlarger{\left|(1-t)y+t\frac{y}{|y|}\right|
-\left|t(\beta_{1}\circ\Phi_{\bar{\rho}})(y,\delta_{1})-t\frac{y}{|y|}\right|}\\
\\
&\geq t+(1-t)|y|-\mathlarger{\frac{t}{4}}
\geq \mathlarger{\frac{1}{2}}+\mathlarger{\frac{t}{4}}>0,\ \forall t\in[0,1].
\end{array}
\end{equation}
By \eqref{eq:Lemma 4.3-(a)(b)} (c) and $\delta_2>1/2$, for $\forall (y,\delta)\in D^{3}$, we have
\begin{equation}\label{eq:D^3-beta_2}
(1-t)\delta_{2}+t(\beta_{2}\circ\Phi_{\bar{\rho}})(y,\delta_{2})>\frac{1}{2},\quad \forall t\in[0,1].
\end{equation}
By \eqref{eq:Lemma 4.3-(a)(b)} (f), for $\forall (y,\delta)\in D^{4}$, we have
\begin{equation}\label{eq:D^4-beta_2}
\Big((1-t)y+t(\beta_{1}\circ\Phi_{\bar{\rho}})(y,\delta),y\Big)_{\mathbb{R}^{3}}>0,\quad \forall t\in[0,1].
\end{equation}
It follows from $\partial D_{y,\delta}=D^{1} \cup D^{2} \cup D^{3}\cup D^{4}$ and \eqref{eq:D^1-beta_2}-\eqref{eq:D^4-beta_2} that \eqref{eq:varphi(y,delta)neq(0,1/2)} holds. So \eqref{eq:d(varphi)=1} is proved.

By \eqref{eq:d(varphi_h)=d(varphi)} and (\ref{eq:d(varphi)=1}), we have
$$
d\left(\varphi_{h},D_{y,\delta},\left(0,0,0,{1}/{2}\right)\right)=1;
$$
that is, the equation $\varphi_{h}(y,\delta)=(0,0,0,{1}/{2})$ has a solution $(\widetilde{y},\widetilde{\delta})\in D_{y,\delta}$. The lemma is proved. $\quad\Box$\\
\textbf{Proof of Theorem 1.1.} Define
\begin{equation}\label{eq:c-define}
c:=\inf_{A\in\Gamma}\sup_{u\in A}I(u),
\end{equation}
$$
K_{c}:=\{u\in \mathcal{N}\cap P:\ I(u)=c,\  I'(u)=0\},
$$
$$
L_{\gamma}:=\{u\in \mathcal{N}:\ I(u)\leq\gamma\}, \quad \gamma\in
\emph{R}.
$$
Choose $\hat{\rho}>0$ small enough such that Lemmas 4.2 and 4.5 hold for $\forall \bar{\rho}<\hat{\rho}$. Now we fix $\bar{\rho}(<\hat{\rho})$. To prove the theorem, it is enough to show that the level $c$
defined by \eqref{eq:c-define} is a critical level, i.e.,
$K_{c}\neq\emptyset$. By Lemma 4.5,  $A\cap
\mathfrak{B}_{0}\neq\emptyset$, $\forall A\in\Gamma$, and by Lemma 4.4, we have
$$
c\geq\inf_{\mathfrak{B}_{0}\cap
P}I(u)=c_{0}>M_{\infty}.
$$
By the choice of $\hat{\rho}$, $\Sigma \in \Gamma$, and Lemma
4.2, we have  $c<2M_{\infty}$. Thus,
$$
M_{\infty}<c<2M_{\infty}.
$$
Suppose by contradiction  $K_{c}=\emptyset$. By Lemma 3.4, it is
easy to see that the Palais-Smale condition holds in
$$
P\cap \mathcal{N}\cap\{u\in D^{1,2}_{0}(\Omega):M_{\infty}<I(u)<2M_{\infty}\}.
$$
Now, using a variant due to Hofer \cite[Lemma 1]{Hofer} of the
classical deformation lemma (see \cite{Rabinowitz}), we can find a
continuous map
$$
\eta:[0,1]\times \mathcal{N}\cap P\rightarrow \mathcal{N}\cap P
$$
and a positive number $\varepsilon_{0}$ such that:\\
$(1)$ $L_{c+\varepsilon_{0}}\backslash L_{c-\varepsilon_{0}}\subset\subset L_{2M_{\infty}}\backslash L_{\frac{c_{0}+M_{\infty}}{2}}$;\\
$(2)$ $\eta(0,u)=u$;\\
$(3)$ $\eta(t,u)=u$, $\forall u\in L_{c-\varepsilon_{0}}\cup\{\mathcal{N}\cap P\backslash L_{c+\varepsilon_{0}}\}$, $\forall t\in[0,1]$;\\
$(4)$ $\eta(1,L_{c+\frac{\varepsilon_{0}}{2}})\subset L_{c-\frac{\varepsilon_{0}}{2}}$.\\
Then let $\widetilde{A}\in \Gamma$ be such that
$$
c\leq\sup_{u\in\widetilde{A}}I(u)<c+\frac{\varepsilon_{0}}{2},
$$
we have $\eta(1,\widetilde{A})\in\Gamma$ and
$$
c\leq\sup_{u\in\eta(1,\widetilde{A})}I(u)<c-\frac{\varepsilon_{0}}{2},
$$
which is a  contradiction.  Thus, $K_c\neq \emptyset$ and  Theorem 1.1 is proved. $\quad\Box$\\


\end{document}